\documentclass[sn-mathphys,Numbered]{sn-jnl}

\usepackage{graphicx}%
\usepackage{multirow}%
\usepackage{amsmath,amssymb,amsfonts}%
\usepackage{amsthm}%
\usepackage{mathrsfs}%
\usepackage[title]{appendix}%
\usepackage{xcolor}%
\usepackage{textcomp}%
\usepackage{manyfoot}%
\usepackage{booktabs}%
\usepackage{algorithm}%
\usepackage{algorithmicx}%
\usepackage{algpseudocode}%
\usepackage{listings}%
\usepackage{comment}%
\usepackage{longtable}%
\usepackage{subfig}
\usepackage{lscape}
\usepackage[shortlabels]{enumitem} 

\usepackage{lscape}

\usepackage{tikz}
\usepackage{pgfplots}
\usetikzlibrary{arrows.meta}
\tikzstyle{block} = [draw, fill=white, rectangle, 
minimum height=2.5em, minimum width=5.em]
\tikzstyle{sum} = [draw, fill=white, circle, node distance=1cm]
\tikzstyle{input} = [coordinate]
\tikzstyle{output} = [coordinate]
\tikzstyle{pinstyle} = [pin edge={to-,thin,black}]
\usetikzlibrary{arrows}
\usetikzlibrary{cd}
\usetikzlibrary{babel}
\usetikzlibrary{matrix}
\usetikzlibrary{arrows.meta}
\tikzstyle{block} = [draw, fill=white, rectangle, 
minimum height=2.5em, minimum width=5.em]
\tikzstyle{sum} = [draw, fill=white, circle, node distance=1cm]
\tikzstyle{input} = [coordinate]
\tikzstyle{output} = [coordinate]
\tikzstyle{pinstyle} = [pin edge={to-,thin,black}]
\usepackage{adjustbox}

\theoremstyle{thmstyleone}%
\newtheorem{theorem}{Theorem}
\theoremstyle{thmstyletwo}%
\newtheorem{example}{Example}%

\theoremstyle{thmstylethree}%
\newtheorem{definition}{Definition}%

\newcommand{\R}{{\mathbb{R}}}
\newcommand{\dd}{{\mathrm{d}}}

\newcommand{\ipopt}{\texttt{IPOPT}}
\newcommand{\snopt}{\texttt{SNOPT}}
\newcommand{\worhp}{\texttt{WORHP}}
\newcommand{\nosnoc}{\texttt{NOSNOC}}
\newcommand{\casadi}{\texttt{CasADi}}
\newcommand{\nosbench}{\texttt{NOSBENCH}}
\newcommand{\nosbenchf}{\texttt{NOSBENCH-F}}
\newcommand{\nosbenchs}{\texttt{NOSBENCH-S}}
\newcommand{\nosbenchrs}{\texttt{NOSBENCH-RS}}
\newcommand{\nosbenchrl}{\texttt{NOSBENCH-RL}}
\newcommand{\nosbenchurl}{\url{https://github.com/apozharski/nosbench}}
\newcommand{\MacMPEC}{\texttt{MacMPEC}}
\newcommand{\MPECLib}{\texttt{MPECLib}}

\newcommand{\NFE}{N_{\mathrm{fe}}}

\newcommand{\Nsucess}{73.8\%}
\newcommand{\NF}{603}
\newcommand{\NS}{100}
\newcommand{\NRS}{32}
\newcommand{\NRL}{167}

\graphicspath{{images}}
\begin{document}

\title[Solving mathematical programs with complementarity constraints arising in nonsmooth optimal control]{Solving mathematical programs with complementarity constraints arising in nonsmooth optimal control}

\author*[1,3]{\fnm{Armin} \sur{Nurkanovi\'c}}\email{armin.nurkanovic@imtek.uni-freiburg.de}

\author[1,3]{\fnm{Anton} \sur{Pozharskiy}}\email{anton.pozharskiy@imtek.uni-freiburg.de}

\author[1,2]{\fnm{Moritz} \sur{Diehl}}\email{moritz.diehl@imtek.uni-freiburg.de}

\affil[1]{\orgdiv{Systems Control and Optimization Laboratory, Department of Microsystems Engineering (IMTEK)}, \orgname{University of Freiburg}}
\affil[2]{\orgdiv{Department of Mathematics}, \orgname{University of Freiburg}}
\affil[3]{\orgdiv{These authors have contributed equally}}

\abstract{
This paper examines solution methods for mathematical programs with complementarity constraints (MPCC) obtained from the time-discretization of optimal control problems (OCPs) subject to nonsmooth dynamical systems.
The MPCC theory and stationarity concepts are reviewed and summarized. 
The focus is on relaxation-based methods for MPCCs, which solve a (finite) sequence of more regular nonlinear programs (NLP), where a regularization/homotopy parameter is driven to zero.
Such methods perform reasonably well on currently available benchmarks.
However, these results do not always generalize to MPCCs obtained from nonsmooth OCPs.
To provide a more complete picture, this paper introduces a novel benchmark collection of such problems, which we call \nosbench.
The problem set includes \NF\ different MPCCs and we split it into a few representative subsets to accelerate the testing.
We compare different relaxation-based methods, NLP solvers, homotopy parameter update and relaxation parameter steering strategies.
Moreover, we check whether the obtained stationary points allow first-order descent directions, which may be the case for some of the weaker MPCC stationarity concepts.
In the best case, the Scholtes' relaxation \cite{Scholtes2001} with \ipopt~\cite{Waechter2006} as NLP solver manages to solve {\Nsucess} of the problems.
This highlights the need for further improvements in algorithms and software for MPCCs. 
}

\maketitle

\section{Introduction}\label{sec:intro}
This paper investigates numerical methods for solving Mathematical Programs with Complementarity Constraints (MPCCs) obtained from the time-discretization of Optimal Control Problems (OCPs) with nonsmooth dynamical systems.
We consider different types of nonsmooth dynamical systems: (a) where the vector field is nonsmooth but continuous, or (b) nonsmooth and discontinuous (e.g., switched systems), and (c) systems with state jumps.
In many cases, the nonsmoothness and combinatorial structure in such systems can be modeled by a coupling of differential algebraic equations with complementarity constraints. 
This gives rise to so-called Dynamic Complementarity Systems (DCSs)~\cite{Brogliato2020}.
A complementarity constraint reads as $0\leq x \perp y \geq0$, which means that 
all entries of the two vectors $x, y \in \R^n$  must be nonnegative, i.e., $x_i \geq0,\ y_i\geq 0$, and that at least one of the components in every pair is zero, i.e., $x_i y_i = 0$. 

A continuous-time optimal control problem (OCP) subject to a DCS has the following form: 
	\begin{subequations}\label{eq:ocp}
	\begin{align}
		\underset{{\substack{x(\cdot), u(\cdot),y(\cdot)}}}{\mathrm{min}} \;   \int_0^T & L(x(t),u(t)) \dd t + E(x(T)) \label{eq:ocp_cost}
		 \\
		\textrm{s.t.} \quad 
		x(0) &= \bar{x}_0\\
		\dot{x}(t) &= f_{\mathrm{c}}(x(t) , y(t), u(t)), \label{eq:ocp_dcs_ode}\\
		0 &= g_{\mathrm{c}}(x(t) , y(t)), \label{eq:ocp_dcs_alg}\\
		0 &\leq G_{\mathrm{c}}(y(t)) \perp H_{\mathrm{c}}(y(t)) \geq 0, \label{eq:ocp_dcs_cc}\\
		0&\leq h_{\mathrm{c}}(x(t) , u(t) ),\; \text{ for almost all } t\in [0,T], \label{eq:ocp_path}\\
		0&\leq r(x(T)),
	\end{align}
\end{subequations}
where $x \in \R^{n_x}$ are the differential states, $y \in \R^{n_y}$ the algebraic states and $u \in \R^{n_u}$ the control inputs.
The function $L: \R^{n_x} \times \R^{n_u} \to \R$ models the stage cost and $E:\R^{n_x}\to \R$ is the terminal cost, $x_0\in\R^{n_x}$ is a given parameter.
The path and terminal constraints are grouped into the functions $h_{\mathrm{c}} : \R^{n_x}  \times \R^{n_u} \to \R^{n_{h_\mathrm{c}}}$ and $r : \R^{n_x}  \to \R^{n_{r}}$, respectively.
The system \eqref{eq:ocp_dcs_ode}-\eqref{eq:ocp_dcs_cc} is a DCS, where
the function $f_{\mathrm{c}}: \R^{n_x} \times \R^{n_y} \to \R^{n_x}$ models the right-hand side of the differential equation, the function 
$g_{\mathrm{c}}: \R^{n_x} \times \R^{n_y} \to \R^{n_{g_{\mathrm{c}}}}$ defines the smooth algebraic equation in the DCS, and the functions $G_{\mathrm{c}}: \R^{n_y} \to \R^{n_{\mathrm{cc}}}$, $H_{\mathrm{c}}: \R^{n_y} \to \R^{n_{\mathrm{cc}}}$ define the complementarity part of the DCS. 
It is assumed that all functions are at least twice continuously differentiable.
Note that even if $L$ and $E$ are smooth and convex, $f_{\mathrm{c}}$, $g_{\mathrm{c}}$ affine, and $-h_{\mathrm{c}}$ and $-r$ smooth convex functions, the complementarity constraints render the OCP \eqref{eq:ocp} nonsmooth and nonconvex.
In principle, it is also possible to have path and terminal constraints on the algebraic variables, which we have omitted to keep the notation light here and below in the discretization.

The DCS \eqref{eq:ocp_dcs_ode}-\eqref{eq:ocp_dcs_cc} abstraction allows one to model a variety of different nonsmooth systems.
Examples are: Filippov differential inclusions reformulated into a DCS via Stewart's \cite{Stewart1990a,Nurkanovic2022} or the Heaviside step reformulation~\cite{Nurkanovic2023,Nurkanovic2023c}, complementarity Lagrangian systems (modeling rigid bodies with friction and impacts)~\cite{Brogliato2020}, relay systems~\cite{Johansson2002}, projected dynamical systems~\cite{Heemels2000a}, Moreau's sweeping processes~\cite{Moreau1977}, and many more~\cite{Brogliato2020}.
Moreover, with the use of time-freezing, several classes of systems with state jumps can be reformulated into Filippov systems~\cite{Nurkanovic2021,Nurkanovic2023a,Nurkanovic2022a,Halm2021}, which in turn lead to DCS.
A detailed overview is given in~\cite{Nurkanovic2023f}.

Here, we {consider} a direct approach, where one first discretizes the continuous-time OCP \eqref{eq:ocp} and obtains a finite-dimensional nonlinear program.
For example, in a direct transcription method, one can use time-stepping methods (e.g., implicit Runge-Kutta (IRK) methods) to discretize the dynamics in time.
In the case of smooth dynamical systems, direct methods are at a very mature stage~\cite{Rawlings2017}. 
However, in the case of nonsmooth systems such as DCS, this approach has some severe limitations. 
In particular, the time-stepping methods, with fixed integration step sizes, have at best first-order accuracy and the derivatives of the solutions/state transition maps with respect to parameters do not converge to the correct values~\cite{Stewart2010,Nurkanovic2020}. 
As a consequence of the wrong derivatives, an optimizer may converge to a spurious solution close to the initial guess~\cite{Stewart2010,Nurkanovic2020}.

These limitations were recently overcome by the Finite Elements with Switch Detection (FESD) method~\cite{Nurkanovic2022,Nurkanovic2023b}, which, inspired by~\cite{Baumrucker2009}, lets the integration step sizes to be degrees of freedom and introduces additional constraints that enable exact switch detection. 
This allows FESD to recover the higher-order accuracy properties of IRK methods and to compute correct sensitivities.
This method is available in the open-source package \nosnoc~\cite{Nurkanovic2022b}.
In this work, we discretize the OCP \eqref{eq:ocp} with the FESD method and obtain a discrete-time OCP.
We introduce a uniform control discretization grid with $N$ intervals and grid points $t_k = t_{k-1} + \frac{T}{N}$. 
Note that the control interval is not the same as the integration interval, and on every control interval $[t_k,t_{k+1}]$ one can apply a FESD method with multiple variable integration steps.
The state approximations are denoted by $x_k \approx x(t_k)$, the control discretization is taken to be constant on the whole interval, i.e., $u(t) = u_k, t\in [t_k,t_{k+1}]$. 

The vector $z_k$ collects all internal variables of the FESD method, e.g., the Runge-Kutta stage variables of the algebraic and differential states,
and the approximations of the algebraic variables $y_k \approx y(t_k)$.
\color{black}
With a slight abuse of notation, they are summarized in the vectors 
$x=\begin{bmatrix} x_0^\top, \ldots,x_N^\top \end{bmatrix}^\top$, 
$z=\begin{bmatrix} z_0^\top, \ldots,z_{N-1}^\top \end{bmatrix}^\top$,
and $u=\begin{bmatrix} u_0^\top, \ldots,u_{N-1}^\top \end{bmatrix}^\top$.
A discretized version of the OCP~\eqref{eq:ocp} reads as~\cite{Nurkanovic2022}:
\begin{subequations}\label{eq:ocp_discrete_time}
	\begin{align}
		\underset{{{x, u, z}}}{\mathrm{min}} \; \sum_{k=0}^{N-1} &\;  \ell(x_k, u_k) + E(x_N) \\
		\textrm{s.t.} \quad 
		x_0 &= \bar{x}_0,\\
		x_{k+1} &= \phi_f(x_k, z_k, u_k ), \; &\text{for all}\ k \in \{0, \ldots, N\!-\!1\},\\
		0 &= \phi_{\mathrm{int}}(x_k, z_k, u_k), \; &\text{for all}\ k \in \{0, \ldots, N\!-\!1\},\\
		0 &\leq \phi_{G}(z_k)  \perp \phi_{H}(z_k) \geq 0, \; &\text{for all}\ k \in \{0, \ldots, N\!-\!1\},\\
		0 &\leq h_{\mathrm{c}}(x_k, u_k),  &\text{for all}\ k \in \{0, \ldots, N\!-\!1\}, \label{eq:ocp_discrete_time_path}\\
		0 &\leq r(x_N).
	\end{align}
\end{subequations}
The functions $\phi_f$, $\phi_{\mathrm{int}}$, $\phi_{G}$ and $\phi_{H}$ define a discrete time DCS are obtained by applying the FESD method to the DCS \eqref{eq:ocp_dcs_ode}-\eqref{eq:ocp_dcs_cc}, for a detailed definition of these functions, cf.~\cite{Nurkanovic2022,Nurkanovic2023c}. 
The terms $\ell(x_k, u_k)$ approximate the stage cost integral in \eqref{eq:ocp_cost} over the intervals $[t_k,t_{k+1}], k = 0,\ldots,N-1$. 
The path constraints \eqref{eq:ocp_path} are for simplicity only evaluated at the points $t_k$, which yields~\eqref{eq:ocp_discrete_time_path}.

By setting $w= \begin{bmatrix} x^\top,z^\top,u^\top \end{bmatrix}^\top \in \R^{n}$ and defining appropriate functions: $f: \R^n \to \R$ for the objective expression,
 $g:\R^n \to \R^{n_g}$ for collecting the equality constraints,  $h:\R^n \to \R^{n_h}$ for the inequality constraints, and $G:\R^n \to \R^m$ and $H:\R^n \to \R^m$ for the complementarity functions, the discrete-time OCP \eqref{eq:ocp_discrete_time} can be compactly written as a generic MPCC:
	\begin{align*}
		\underset{w\in \R^{n}}{\mathrm{min}} \;  \quad &f(w)\\
		\textrm{s.t.} \quad 
		&g(w)=0,\\
		&h(w)\geq0, \\
		&0 \leq G(w) \perp H(w) \geq 0.
	\end{align*}

 Mathematical Programs with Equilibrium Constraints (MPEC) are NLPs that have a parametric variational inequality or optimization problem as an  constraint~\cite{Luo1996}. 
 Such constraints can be under suitable conditions replaced by equivalent complementarity conditions. 
 However, in the literature, because of the easier pronunciation, the acronym MPEC is frequently used for the problem above.
\color{black}
There are few equivalent ways to state the complementarity constraints  $0 \leq G(w) \perp H(w) \geq 0$
\color{black} as formally smooth constraints:
\begin{enumerate}
	\item $ G(w)\geq 0,\ H(w) \geq 0,\ G_i(w)H_i(w) \leq 0,\ \textit{for all} \ i\in\{1,\ldots m\}$,
	\item $ G(w)\geq 0,\  H(w) \geq 0,\  G(w)^\top H(w) \leq 0$,
	\item $ G(w)\geq 0,\ H(w) \geq 0,\  G(w)^\top H(w) = 0$,
	\item $ G(w)\geq 0,\ H(w) \geq 0,\  G_i(w)H_i(w) = 0,\ \textit{for all} \ i\in\{1,\ldots m\}$,
	\item $\Phi_{\mathrm{C}}(G(w),H(w))=0$.
\end{enumerate}
In (4), $\Phi_{\mathrm{C}}$ is a so-called C-function~\cite{Facchinei2003}, which has the property $\Phi_{\mathrm{C}}(G(w),H(w))=0$ if and only if $0\leq G(w) \perp H(w) \geq 0$. 
C-functions can be smooth or nonsmooth.
To be consistent with most of the MPCC literature, we will work with complementarity constraints written via the inequality constraints in (1):
\begin{subequations}\label{eq:mpcc}
	\begin{align}
		\underset{w\in \R^{n}}{\mathrm{min}} \;  \quad &f(w)\\
		\textrm{s.t.} \quad 
		&g(w)=0,\\
		&h(w)\geq0, \label{eq:mpcc_ineq}\\
		&G_i(w)\geq 0,\ H_i(w) \geq 0,\ G_i(w) H_i(w)  \leq 0, &\textit{for all} \ i\in\{1,\ldots m\}.  \label{eq:mpcc_cc}
	\end{align}
\end{subequations}
\color{black}
There are no significant theoretical differences or computational advantages in using one of the other equivalent forms.

It is very common to introduce slack variables for the functions $G(w)$ and $H(w)$ to have only linear functions in the complementarity conditions, i.e., instead of \eqref{eq:mpcc_cc}, one has:
\begin{align*}
	s_G = G(w), \; s_H = H(w),\ s_G\geq 0,\; s_H \geq 0,\ s_{G,i} s_{H,i} \leq 0,\ \ \text{for all}\ i \in \{1,\ldots,m\}.
\end{align*}
This is called the \textit{vertical form} of the MPCC.
This does not change any of the theoretical considerations and we stick to the notation of most of the MPCC literature where the slacks are not introduced.
However, for the efficacy of numerical solvers, it is often beneficial to introduce the slacks~\cite{Fletcher2006}, and we do so in the numerical experiments.

At this point, we mention also the class of Mathematical Programs with Vanishing Constraints (MPVC), where the complementarity constraints \eqref{eq:mpcc_cc} are replaced by $G_i(w) H_i(w)\geq0, H_i(w)\geq 0,\ \text{for all}\ i = \{1,\ldots,m\}$~\cite{Achtziger2008}.
They can be reformulated into equivalent MPCC, but it is often numerically more beneficial to treat them directly~\cite{Achtziger2008}.
MPVCs arise often in the relaxation of OCPs with integer control and combinatorial constraints~\cite{Kirches2010y,Jung2013}. 

The MPCC \eqref{eq:mpcc} is a nonlinear program (NLP) for which we need efficient and robust numerical solution methods.
If a point $w^*$ satisfies a Constraint Qualification (CQ), e.g., the Linear Independence Constraint Qualification (LICQ), then the Karush-Kuhn-Tucker (KKT) conditions are necessary for $w^*$ to be local minimizer of~\eqref{eq:mpcc}~\cite{Nocedal2006}.
Standard nonlinear programming algorithms solve the KKT conditions to find a solution candidate.
Unfortunately, due to the complementarity constraints \eqref{eq:mpcc_cc}, standard CQs, such as LICQ and Mangasarian-Fromovitz Constraint Qualification (MFCQ) are violated at all feasible points~\cite{Scheel2000,Ye2005}.
This implies both numerical and theoretical difficulties. 
On the one hand, the violation of MFCQ means that the set of Lagrange multipliers is necessarily unbounded, which leads to computational difficulties~\cite{Izmailov2014}.
On the other hand, because of the absence of CQs, the KKT conditions may no longer be necessary for optimality anymore. 
This has led to the development of a tailored theory and solution methods for MPCCs, which we recall in detail in the following two sections.

A good way to assess the performance and robustness of a tailored MPCC method is to apply it to a benchmark set.
Two widely used benchmarks for MPCCs are the \MPECLib~\cite{Dirkse2004} and \MacMPEC~\cite{Leyffer2000}.
Baumrucker et al.~\cite{Baumrucker2008} test several MPCC methods on the \MPECLib\ test set, and examples from optimal process control. 
Hoheisel et al. \cite{Hoheisel2013} compare several relaxation-based methods on the \MacMPEC\ problem set.
All these experiments report success rates above 90\%. 
In MPCCs arising from nonsmooth OCPs, we did not observe such robustness and high success rates, which motivated us to introduce a new collection of test problems to assess the performance of MPCC methods.

\paragraph{Contributions}
To learn more about the performance of MPCC methods in solving discrete-time nonsmooth OCPs, we introduce the \nosbench~problem set.
It contains {\NF} MPCCs generated via~\nosnoc~\cite{Nurkanovic2022b} from {33} continuous-time OCP and simulation problems. 
Together with a review of the theory and MPCC solution methods, the introduction of \nosbench\ is the main contribution of this paper.
Furthermore, we compare {nine} different relaxation-based methods from the literature together with three NLP solvers: \ipopt~\cite{Waechter2006}, \snopt~\cite{Gill2005} and \worhp~\cite{Bueskens2013}. 
We also compare homotopy parameter update and steering strategies.
Some of the weaker multiplier-based MPCC stationarity concepts may allow first-order descent directions.
We check is this the case for the solutions computed in our experiments.
In our experiments, the Scholtes' relaxation \cite{Scholtes2002} with \ipopt~\cite{Waechter2006} as NLP subproblem solver is the most successful method-solver combination and solves {\Nsucess}\ of the problems in the full \nosbench\ problem set.
Furthermore, we validate the correctness of our implementations by running them on the {\MacMPEC} test set, where we obtain results that aligns with those reported in literature~\cite{Thierry2020,Izmailov2012b,Fletcher2002a}.

\paragraph{Outline}
Section~\ref{sec:mpcc_theory} reviews briefly the standard MPCC theory. 
In Section~\ref{sec:mpcc_methods}, we review easy-to-implement MPCC solution methods and comment on some other promising methods, which yet lack good implementations.
Section \ref{sec:nosbench} introduces the \nosbench\ problem test set and Section~\ref{sec:results} discusses the results we obtain.
We summarize our findings in Section~\ref{sec:conclusions}.

\section{Optimality conditions for MPCCs}\label{sec:mpcc_theory}
Due to the violation of the CQs, the standard NLP theory is often not applicable to MPCCs in the form of Eq. \eqref{eq:mpcc}.
This has several negative consequences: 
(a) the set of Lagrange multipliers is necessarily unbounded,
(b) the gradients of the active constraints are linearly dependent at all feasible points, and 
(c) the linearization of \eqref{eq:mpcc} can be inconsistent arbitrarily close to a stationarity point~\cite{Fletcher2005}.

We review first-order necessary optimality conditions and several stationarity concepts for MPCCs. 
Some {optimality} conditions are purely geometric, others rely on Lagrange multipliers. 
If there exists an $i$ such that $G_i(w) =0, H_i(w) = 0$, the multiplier-based stationary concepts may not be strong enough to characterize local minimizers.
The {MPCC-tailored theory} presented here has its origins in~\cite{Flegel2005a,Fletcher2006,Ye2005,Luo1996,Scheel2000}.

The feasible set of the MPCC \eqref{eq:mpcc} is denoted by $\Omega_{\mathrm{MPCC}} = \{ w \in \R^{n} \mid 
g(w) = 0,\ h(w) \geq0,\ G(w) \geq 0,\ H(w) \geq 0,\ G_i(w)H_i(w) \leq 0,  \ \text{for all}\ i \in \{1,\ldots,m\}$.
We define the following index sets which depend on a feasible point $w \in \Omega_{\mathrm{MPCC}}$:
\begin{align*}
	\mathcal{I}_{+0}(w) &=	\{i \in \{1,\ldots,m\} \mid G_i(w)>0, H_i(w)=0\},\\
	\mathcal{I}_{0+}(w) &= \{i \in \{1,\ldots,m\}\mid G_i(w)=0, H_i(w)>0\},\\
	\mathcal{I}_{00}(w) &= \{i \in \{1,\ldots,m\}\mid G_i(w)=0, H_i(w)=0\}.
\end{align*}
For ease of notation, if clear from the context, we omit the argument in the index sets.
The set $\mathcal{I}_{00}$ is called the set of degenerate indices and is the source of most theoretical and numerical difficulties. 
If $\mathcal{I}_{00}$ is empty, we say that a solution $w^*$ satisfies \textit{strict complementarity}. 
This notion should not be confused with the notion of strict complementarity of an inequality constraint \eqref{eq:mpcc_ineq} and the corresponding Lagrange multiplier in the NLP.

For a closed set $\Omega$ and a point $w \in \Omega$, a vector $d$ is said to be tangent to $\Omega$ at $w$ if there exists a sequence $w_k \in \Omega$ with $w_k \to w$, along with a sequence of positive scalar $t_k \to 0$ such that $d_k = \frac{w_k-w}{t_k} \to d$. 
The set of all vectors $d$ that are tangent to $\Omega$ at a point $w \in \Omega$ is called the tangent cone to $\Omega$ at $w$ and denoted by $\mathcal{T}_{\Omega}(w)$. 
Hence, we denote by $\mathcal{T}_{\Omega_{\mathrm{MPCC}}}(w)$ the tangent cone of a feasible point of \eqref{eq:mpcc}. 
The active set of the inequality constraints $h(w)\geq0$ is defined as  the set $\mathcal{A}(w) = \{i \in \{1,\ldots,n_h\} \mid h_i(w)=0\}$.
The linearized feasible cone $\mathcal{F}_{\Omega_{\mathrm{MPCC}}}(w)$ of the MPCC \eqref{eq:mpcc} at a feasible point $w \in \Omega_{\mathrm{MPCC}}$ is defined as the set:
\begin{align*}
	\begin{split}
		\mathcal{F}_{\Omega_{\mathrm{MPCC}}}(w) = \{ d \in \R^n \mid 
		&\nabla g(w)^\top d =0,\\
		&\nabla h_{i} (w)^\top d \geq 0, 	\textrm{ for all } i \in \mathcal{A}(w),\\
		&\nabla G_{i} (w)^\top d =  0, 		\textrm{ for all } i \in \mathcal{I}_{0+}(w),\\ 
		&\nabla H_{i} (w)^\top d =  0, 		\textrm{ for all } i \in \mathcal{I}_{+0}(w),\\ 
		&\nabla G_{i} (w)^\top d \geq  0, 	\textrm{ for all } i \in \mathcal{I}_{00}(w),\\ 
		&\nabla H_{i} (w)^\top d \geq  0, 	\textrm{ for all } i \in \mathcal{I}_{00}(w)\}.
	\end{split}
\end{align*}
This set is a polyhedral convex cone.
On the other hand, if $\mathcal{I}_{00} \neq \emptyset$, then the tangent cone $\mathcal{T}_{\Omega_{\mathrm{MPCC}}}(w)$ is, due to the complementarity constraints, a union of polyhedral cones. 
Consequently, it is possibly a nonconvex cone. 
To see this, regard the set $ \Omega = \{ (w_1,w_2) \in \R^2 \mid {w_1\geq0},\ {w_2 \geq 0},\ w_1w_2\leq0\}$. 
It can be verified that $\mathcal{T}_{\Omega}(0) = \Omega$, which is nonconvex, and $\mathcal{F}_{\Omega}(0) = \{ (d_1,d_2) \in \R^2 \mid d_1\geq0,\ d_2 \geq 0\}$, which is convex.
Therefore, the linearized feasible cone cannot locally capture the structural nonconvexity of the complementarity constraints.
The KKT conditions require a CQ to hold in order to be necessary for optimality.
We can see that even the rather nonrestrictive Abadie CQ (ACQ), which simply requires $\mathcal{F}_{\Omega_{\mathrm{MPCC}}}(w) = \mathcal{T}_{\Omega_{\mathrm{MPCC}}}(w)$, cannot be expected to hold for MPCCs~\cite{Flegel2005}.
Only the weaker Guignard CQ (GCQ)~\cite{Guignard1969}, which requires that the polar cones of $\mathcal{F}_{\Omega_{\mathrm{MPCC}}}(w)$ and $\mathcal{T}_{\Omega_{\mathrm{MPCC}}}(w)$ are equal, has a chance to hold~\cite{Flegel2006}. 
However, it is difficult to verify this condition in practice. 

To have a more powerful theoretical tool, the MPCC linearized feasible cone $\mathcal{F}_{\Omega_{\mathrm{MPCC}}}^{\mathrm{MPCC}}(w)$ can be used~\cite{Flegel2005a,Pang1999,Scheel2000}.
This cone is defined at a feasible point $w$ as
	\begin{align*}
		\begin{split}
			\mathcal{F}_{\Omega_{\mathrm{MPCC}}}^{\mathrm{MPCC}}(w) = \{ d \in \R^n \mid 
			&\nabla g(w)^\top d =0,\\
			&\nabla h_{i} (w)^\top d \geq 0, \textrm{ for all } i \in \mathcal{A}(w),\\
			&\nabla G_{i} (w)^\top d = 0, \textrm{ for all } i \in \mathcal{I}_{0+}(w),\\
			&\nabla H_{i} (w)^\top d = 0, \textrm{ for all } i \in \mathcal{I}_{+0}(w),\\
			&0 \leq  \nabla G_{i} (w)^\top d \perp \nabla H_{i} (w)^\top d \geq0, \textrm{ for all } i \in \mathcal{I}_{00}(w)
			\}.
		\end{split}
	\end{align*}
The combinatorial structure is kept for the degenerate index set $\mathcal{I}_{00}$ and the cone 	$\mathcal{F}_{\Omega_{\mathrm{MPCC}}}^{\mathrm{MPCC}}(w)$ is nonconvex, if $\mathcal{I}_{00}$ is nonempty.
In example from above, it holds that $\mathcal{F}^{\mathrm{MPCC}}_{\Omega}(0) = \mathcal{T}_{\Omega}(0)$.

We proceed with stating optimality conditions and defining stationarity concepts for MPCCs.
First-order necessary optimality conditions can be stated in terms of the tangent cone.
\begin{theorem}
Let $w^*\in\Omega_{\mathrm{MPCC}}$ be a local minimizer of \eqref{eq:mpcc}, then it holds that
\begin{align}\label{eq:geometric_b_stationarity}
	\nabla f(w^*)^\top d \geq 0\; \textrm{ for all } d \in \mathcal{T}_{\Omega_{\mathrm{MPCC}}}(w^*).
\end{align}
\end{theorem}
If a point $w^*$ satisfies the condition above, in the MPCC literature it is said that geometric Bouligand stationarity (geometric B-stationarity) holds~\cite{Flegel2005a,Luo1996}.
For computational purposes, algebraic stationarity concepts are more useful.
The algebraic Bouligand stationarity (or just B-stationarity) \cite{Luo1996,Scheel2000} reads as follows.
\begin{definition}[B-stationarity]\label{def:b_stationariry}
 A feasible point $w^*\in\Omega_{\mathrm{MPCC}}$ of the MPCC \eqref{eq:mpcc} satisfies B-stationarity if it holds that
\begin{align}\label{eq:algebraic_b_stationarity}
	\nabla f(w^*)^\top d \geq 0\; \textrm{ for all } d \in \mathcal{F}_{\Omega_{\mathrm{MPCC}}}^{\mathrm{MPCC}}(w^*).
\end{align}
\end{definition} 
Or equivalently \cite{Luo1996,Scheel2000}, a point $w^*\in\Omega_{\mathrm{MPCC}}$ is called B-stationary if $d = 0$ is a local minimizer of the following linear program with complementarity constraints:
\begin{align}\label{eq:b_stationariry}
	\min_{d \in \R^{n}} \quad &   \nabla f(w^*)^\top d\quad 	\mathrm{s.t.} \quad   d \in \mathcal{F}_{\Omega_{\mathrm{MPCC}}}^{\mathrm{MPCC}}(w^*).
\end{align}
It was shown in \cite{Flegel2005a} that $\mathcal{T}_{\Omega_{\mathrm{MPCC}}}(w) \subseteq \mathcal{F}_{\Omega_{\mathrm{MPCC}}}^{\mathrm{MPCC}}(w^*)$.
This means that B-stationarity is less restrictive and it implies geometric B-stationarity.
However, the converse is not always true, and we discuss below conditions when this is the case and when B-stationarity is necessary for optimality.
B-stationarity is expensive to verify, as it requires the solution of a nonconvex optimization problem.
In the worst case, this may require the solution of an exponential number of linear programs, unless some stronger regularity conditions hold.

As in the standard NLP theory, we may want to use Lagrange multiplier-based stationarity concepts, where we hopefully do not need to solve a combinatorial problem to find a solution candidate.
As the standard CQs do not hold, we cannot use the KKT conditions of the MPCC \eqref{eq:mpcc}, but we apply them to some auxiliary Nonlinear Programs (NLPs). 
The stationarity conditions for MPCCs are derived from more regular NLPs (defined next) which are locally associated with the initial MPCC \eqref{eq:mpcc}~\cite{Scheel2000}.
\begin{definition}[Auxiliary NLP] \label{def:mpcc_nlps}
	Let ${w}^* \in \Omega_{\mathrm{MPCC}}$. We define the following auxiliary NLPs:
	\begin{itemize}
		\item {the Relaxed NLP (RNLP) for ${w}^*\in \Omega_{\mathrm{MPCC}}$ is defined as}
		\allowdisplaybreaks
		\begin{subequations} \label{eq:rnlp}
			\begin{align}
				\min_{x \in \R^{n}} \quad &   f(w) \label{eq:rnlp_objective}\\
				\mathrm{s.t.} \quad  &g(w) = 0 \label{eq:rnlp_eq},\\
				&h(w)\geq 0, \label{eq:rnlp_ineq}\\
				&G_i(w) = 0,\ H_i(w) \geq 0,\; i \in 	\mathcal{I}_{0+}(w^*), \label{eq:rnlp_branch_h}\\
				&G_i(w) \geq 0,\ H_i(w) = 0,\; i \in 	\mathcal{I}_{+0}(w^*), \label{eq:rnlp_branch_g}\\
				&G_i(w) \geq 0,\ H_i(w) \geq 0,\; i \in 	\mathcal{I}_{00}(w^*) \label{eq:rnlp_biactive},
			\end{align}
		\end{subequations}
		\item {the Tight NLP (TNLP) for ${w}^*\in \Omega_{\mathrm{MPCC}}$ is defined as} the RNLP \eqref{eq:rnlp}, except that the constraints \eqref{eq:rnlp_biactive} are replaced by:
		\begin{align*}
			&G_i(w) = 0,\ H_i(w) = 0,\; i \in 	\mathcal{I}_{00}(w^*).
		\end{align*}
		\item Let $(\mathcal{I}_1,\mathcal{I}_2)$ be a partition of $\mathcal{I}_{00}$ such that $\mathcal{I}_1\cup\mathcal{I}_2= \mathcal{I}_{00}$ and $\mathcal{I}_1\cap\mathcal{I}_2 = \emptyset$. 
		The Branch NLP ($\mathrm{BNLP}_{(\mathcal{I}_1,\mathcal{I}_2)}$) for ${w}^*\in \Omega_{\mathrm{MPCC}}$ is 
		defined as the RNLP \eqref{eq:rnlp}, except that the constraints \eqref{eq:rnlp_branch_h}-\eqref{eq:rnlp_biactive} are now replaced by:
		\begin{align*}
		&G_i(w) \geq 0,\ H_i(w) = 0,\; i \in 	\mathcal{I}_{+0}(w^*) \cup \mathcal{I}_1(w^*), \\
		&G_i(w) = 0,\ H_i(w) \geq 0,\; i \in 	\mathcal{I}_{0+}(w^*)\cup \mathcal{I}_2(w^*).
		\end{align*}
	\end{itemize}
\end{definition}
\begin{figure}[t]
	\centering
	{\includegraphics[width=\textwidth]{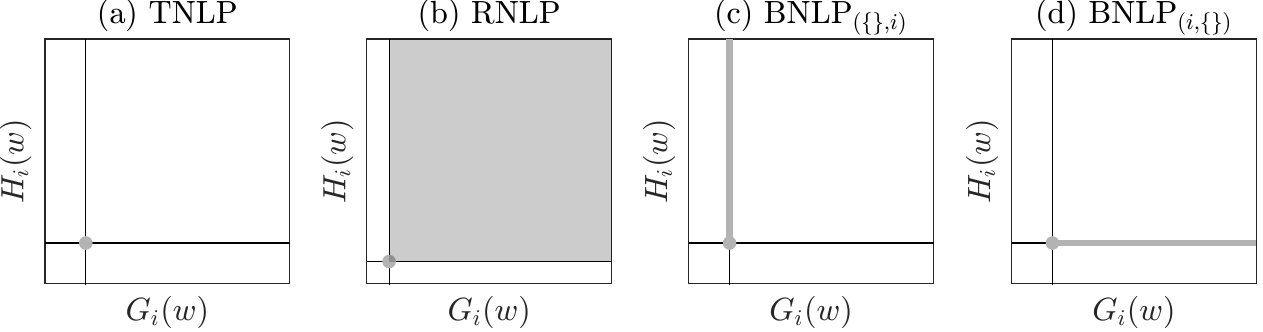}}
	\caption{Feasible sets of the auxiliary NLPs as defined in Def.~\ref{def:mpcc_nlps}.}
	\label{fig:mpcc_nlp}
\end{figure}
Figure~\ref{fig:mpcc_nlp} illustrates the feasible sets of the auxiliary NLPs for $i \in \mathcal{I}_{00}$.
We denote the feasible sets of the RNLP, TNLP, and $\mathrm{BNLP}_{(\mathcal{I}_1,\mathcal{I}_2)}$ by $\Omega_{\mathrm{RNLP}}$, $\Omega_{\mathrm{TNLP}}$ and $\Omega_{\mathrm{BNLP}_{(\mathcal{I}_1,\mathcal{I}_2)}}$, respectively.

The usual NLP concepts such as first-order optimality conditions, stationary points, second-order conditions, and constraint qualification for MPCC are defined in terms of these auxiliary NLP. 
To see that this approach makes sense we look at how these problems and their solutions are related. 
It is not difficult to see that the following holds for $w^* \in \Omega_{\mathrm{MPCC}}$ \cite{Scheel2000}:
\begin{align}\label{eq:mpcc_feasible_sets}
	\Omega_{\mathrm{TNLP}}  = \bigcap_{(\mathcal{I}_1,\mathcal{I}_2)} \Omega_{\mathrm{BNLP}_{(\mathcal{I}_1,\mathcal{I}_2)}} \subset \Omega_{\mathrm{MPCC}}
	= \bigcup_{(\mathcal{I}_1,\mathcal{I}_2)} \Omega_{\mathrm{BNLP}_{(\mathcal{I}_1,\mathcal{I}_2)}}  \subseteq \Omega_{\mathrm{RNLP}}.
\end{align}
The same relations hold for the corresponding tangent cones at $w^*$ as well.
Furthermore, for a feasible point of the MPCC \eqref{eq:mpcc} $w^* \in \Omega_{\mathrm{MPCC}}$ the following can be said \cite{Scheel2000}. 
If $w^*$ is a local minimizer of the RNLP, then it is a local minimizer of the MPCC. 
The converse is not true. 
If $w^*$ is a local minimizer of the MPCC then it is a local minimizer of the TNLP. 
The point $w^*$ is a local minimizer of the MPCC if and only if it is a local minimizer of every $\mathrm{BNLP}_{(\mathcal{I}_1,\mathcal{I}_2)}$. 
The last assertion once again highlights the combinatorial nature of MPCCs, since $2^{|\mathcal{I}_{00}|}$ branch NLPs must be checked to make conclusions about optimality.
Fortunately, as we will see below, under reasonable assumptions we do not have to check every branch NLP but only the RNLP or TNLP to characterize a stationary point of the MPCC.

All these difficulties arise because of the degenerate indices $i \in \mathcal{I}_{00}$. 
If this set is empty, all auxiliary NLPs collapse to the same problem, and there is no combinatorial structure due to the BNLP anymore.
It can be seen that the tangent cone of $\Omega_{\mathrm{MPCC}}$ will be convex since there will be no rays that start from the degenerate point.
Assuming that other constraints in the MPCC do not cause the violation of the ACQ, then we have that the standard ACQ holds for the MPCC and thus we can apply the KKT-conditions to verify the stationarity of $w^*\in \Omega_{\mathrm{MPCC}}$ in this fortunate case. 
In other words, $w^* \in \Omega_{\mathrm{MPCC}}$ is a local minimizer of the MPCC if and only if it is a local minimizer of the RNLP/TNLP, which are equal in this case~\cite{Scheel2000}.

Next, we define the MPCC-specific Lagrangian, CQs, and stationarity concepts. 
The MPCC Lagrangian is the \textit{standard} Lagrangian for the RNLP/TNLP, and reads as:
	\begin{align}
		\mathcal{L}^{\mathrm{MPCC}}(x,\lambda,\mu,\nu,\xi) = f(w) - \lambda^\top g(w) - \mu^\top h(w)- \nu^\top G(w)- \xi^\top H(w),
	\end{align}
	with the MPCC Lagrange multipliers $\lambda \in \R^{n_g}$, $\mu \in \R^{n_h}$, $\nu \in \R^{m}$ and $\xi \in \R^{m}$.
It differs from the standard Lagrangian for the MPCC \eqref{eq:mpcc} in omitting the bilinear terms $G_i(w)H_i(w)\leq0$ and their multipliers.

Next we define some tailored MPCC CQs.
\begin{definition}
	The MPCC \eqref{eq:mpcc} is said to satisfy the MPCC-LICQ (MPCC-MFCQ) at a feasible point $w^*$ if the corresponding $\mathrm{TNLP}$ for $w^*$ satisfies the LICQ (MFCQ) at the same point $w^*$.
\end{definition}
The linearized feasible cone of a TNLP is always convex, and as seen in the discussion at the beginning of this section we can expect the standard ACQ to be violated if $\mathcal{I}_{00} \neq \emptyset$.
This motivated the definition of the MPCC-ACQ and MPCC-Guignard CQ (MPCC-GCQ) in terms of the nonconvex cones~\cite{Flegel2005a,Flegel2005}.
First, recall that given a cone $\mathcal{K} \subseteq \R^n$, its polar cone is defines as $\mathcal{K}^\circ = \{ d \in \R^n \mid  d^\top v \leq 0, \textrm{ for all } v \in \mathcal{K} \}$.
\begin{definition}
The MPCC-ACQ (MPCC-GCQ) holds at $w^* \in \Omega_{\mathrm{MPCC}}$ if and only if 
$\mathcal{F}_{\Omega_{\mathrm{MPCC}}}^{\mathrm{MPCC}}(w^*)= \mathcal{T}_{\Omega_{\mathrm{MPCC}}}(w^*)$ ($\mathcal{F}_{\Omega_{\mathrm{MPCC}}}^{\mathrm{MPCC}}(w^*)^\circ= \mathcal{T}_{\Omega_{\mathrm{MPCC}}}(w^*)^\circ$).
\color{black}
\end{definition}
Similar as for the standard CQ, for the MPCC CQs the following implications hold~\cite{Ye2005,Schwartz2011}: 
${\textrm{MPCC-LICQ} \implies \textrm{MPCC-MFCQ} \implies \textrm{MPCC-ACQ} \implies \textrm{MPCC-GCQ}}.$ 

Inspired by the KKT conditions for standard non-degenerate NLPs, several stationarity concepts that rely on the auxiliary NLPs and their Lagrange multipliers can be defined~\cite{Scheel2000}. 
If an appropriate MPCC-CQ holds, they are necessary for optimality, as we will discuss below.
\begin{definition}[Stationarity concepts for MPCCs]\label{def:mpcc_stationarity}
	Let $w^*$ be feasible for the MPCC~\eqref{eq:mpcc}.
	\begin{itemize}
		\item Weak Stationarity (W-stationarity) \cite{Scheel2000}: A point $w^*\in\Omega_{\mathrm{MPCC}}$ is called W-stationary if the corresponding TNLP admits the satisfaction of the KKT conditions, i.e., there exist Lagrange multipliers $\lambda^*,\mu^*,\nu^*$ and $\xi^*$ \color{black} such that:
		\begin{align*}
			&\nabla_x \mathcal{L}^{\mathrm{MPCC}}(w^*,\lambda^*,\mu^*,\nu^*,\xi^*) = 0,\\
			&g(w^*) = 0,\\
			&0 \leq \mu^* \perp h(w^*) \geq 0,\\
			& G(w^*) \geq0, \nu^*_i = 0, \; \text{ for all } i \in 	\mathcal{I}_{+0}(w^*),\\
			& H(w^*) \geq0, \xi^*_i = 0,\; \text{ for all } i \in 	\mathcal{I}_{0+}(w^*),\\
			& G_i(w^*) =0,\; \nu_i^* \in \R, \; \text{ for all } i \in 	\mathcal{I}_{0+}(w^*)\cup\mathcal{I}_{00}(w^*),\\
			& H_i(w^*) =0,\; \xi_i^* \in \R, \; \text{ for all } i \in 	\mathcal{I}_{+0}(w^*)\cup\mathcal{I}_{00}(w^*).
		\end{align*}
		\item Strong Stationarity (S-stationarity) \cite{Scheel2000}:  A point $w^*\in\Omega_{\mathrm{MPCC}}$ is called S-stationary if it is weakly stationary and $\nu^*_i \geq 0, \xi^*_i \geq0$ for all $i \in \mathcal{I}_{00}(w^*)$. 
		In other words, it is a KKT point of the corresponding RNLP.
		\item Clarke Stationarity (C-stationarity) \cite{Scheel2000}:  A point $w^*\in\Omega_{\mathrm{MPCC}}$ is called C-stationary if it is weakly stationary and $\nu^*_i\xi^*_i \geq0$ for all $i \in \mathcal{I}_{00}(w^*)$.
		\item Mordukhovich Stationarity (M-stationarity) \cite{Scheel2000}:  A point $w^*\in\Omega_{\mathrm{MPCC}}$ is called M-stationary if it is weakly stationary and if either $\nu^*_i >0$ and $\xi^*_i >0$ or $\nu^*_i\xi^*_i =0$ for all $i \in \mathcal{I}_{00}(w^*)$.
		\item Abadie Stationarity (A-stationarity) \cite{Flegel2005a}: A point $w^*\in\Omega_{\mathrm{MPCC}}$ is called A-stationary if it is weakly stationary and $\nu^*_i \geq 0$ or $\xi^*_i \geq0$ for all $i \in \mathcal{I}_{00}(w^*)$.
	\end{itemize}
\end{definition}
\begin{figure}[t]
	\centering
	{\includegraphics[width=\textwidth]{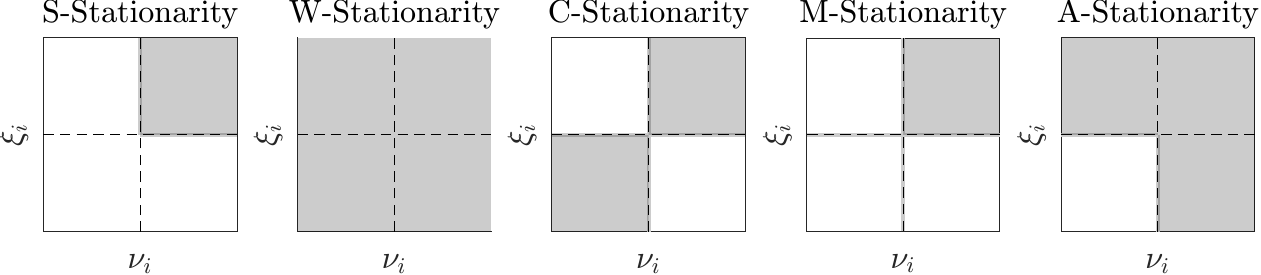}}
	\vspace{-0.5cm}
	\caption{Sign restrictions for MPCC multiplier in different stationarity concepts.}
	\label{fig:mpcc_stationarirty}
\end{figure}
The feasible sets for the MPCC multipliers $\nu^*$ and $\xi^*$ are depicted in Figure \ref{fig:mpcc_stationarirty}. 
Observe that if $\mathcal{I}_{00} = \emptyset$, then all multiplier-based stationarity concepts collapse to the same. 

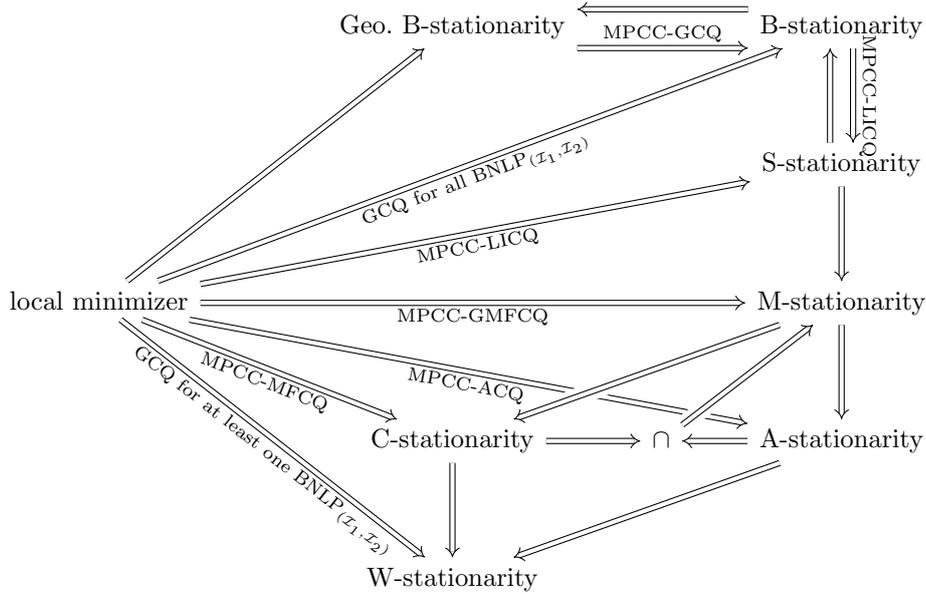
\begin{figure}
	\begin{tikzcd}
		&  & \textrm{Geo. B-stationarity}  \arrow[rr,Leftarrow,shift left=1.5ex] 
		\arrow[rr,"\textrm{MPCC-GCQ}",Rightarrow,shift right=1.9ex] 		
		&	 & \textrm{B-stationarity}  
		\arrow[dd,Rightarrow,"\textrm{MPCC-LICQ}",above,sloped,shift left=1.0ex] 
		\arrow[dd,Leftarrow,,shift right=1.0ex] 
		\\
		&   &   &  & \\ 
		&   &  &   & \textrm{S-stationarity}	
		\arrow[dd,Rightarrow] 
		\\
		&   &   &  & \\ 
		\textrm{local minimizer} 
		\arrow[uuuurr,Rightarrow] 
		\arrow[uuuurrrr,Rightarrow,"\textrm{GCQ for all BNLP$_{(\mathcal{I}_1,\mathcal{I}_2)}$}",sloped,swap] 
		\arrow[uurrrr,Rightarrow,"\textrm{MPCC-LICQ}",sloped,swap] 
		\arrow[rrrr,Rightarrow,"\textrm{MPCC-GMFCQ}",sloped,swap] 
		\arrow[rrdd,Rightarrow,Rightarrow,"\textrm{MPCC-MFCQ}",sloped,swap] 
		\arrow[rrrrdd,Rightarrow,"\textrm{MPCC-ACQ}",sloped,swap] 
		\arrow[rrdddd,Rightarrow,"\textrm{GCQ for at least one BNLP$_{(\mathcal{I}_1,\mathcal{I}_2)}$}",sloped,swap] 
		&   &  &   & \textrm{M-stationarity}	
		\arrow[dd,Rightarrow] 
		\arrow[ddll,Rightarrow,crossing over] 
		\\
		&   &   &  & \\ 
		&   &\textrm{C-stationarity} \arrow[dd,Rightarrow] \arrow[r,Rightarrow] & \cap \arrow[r,Leftarrow]  \arrow[ruu, Rightarrow, crossing over] & \textrm{A-stationarity}	\arrow[lldd,Rightarrow] \\ 
		&   &   &  & \\ 
		&   &      \textrm{W-stationarity} & & \\
	\end{tikzcd}
	\caption{Diagram summarizing MPCC-CQs and necessity of different stationarity concepts for optimality.
	This diagram was inspired by \cite [Theorem 5.13]{Lenders2018} and \cite[Figure 2]{Kim2020}. }
	\label{fig:mpcc_relations}
\end{figure}

The many different stationarity concepts might be confusing and one might be wondering if some of them are needed at all.
However, these stationarity concepts are crucial for studying numerical methods for MPCC. 
As we will see in the next section, MPCCs are usually solved by solving a (finite) sequence of related and more regular NLPs. 
Depending on the underlying assumptions, the accumulation points of these methods are some of the stationary points defined above. 
Therefore, it is important to understand under which conditions these stationarity concepts are indeed necessary for optimality. 
It turns out that all of them can be necessary for local optimality if some additional specialized CQs hold~\cite{Fletcher2006,Scheel2000}. 
The results from \cite{Fletcher2006,Ye2005,Luo1996,Outrata1999,Scheel2000} are summarized in the diagram in Figure \ref{fig:mpcc_relations}.
For the missing CQ definitions see \cite{Ye2005,Schwartz2011}.


We make a few comments on the relations in Figure \ref{fig:mpcc_relations}.
As already discussed, if the ACQ does not hold, the KKT conditions are not applicable for characterizing a geometric B-stationary point. 
B-stationarity always implies geometric B-stationarity, the converse requires additionally the MPCC-GCQ to hold \cite{Flegel2005a}. 
Note that by Definition \ref{def:b_stationariry}, B-stationarity does not allow a first-order descent direction.
Weaker concepts, mainly the multiplier-based ones, may allow first-order descent directions, even if the MPCC-LICQ holds. 
We illustrate this in the next example from \cite{Scheel2000}.
\begin{example}[Descent directions for C-stationarity]
	Regard the MPCC: 
	\begin{align*}
		\min_{w_1,w_2 \in \R} (w_1-1)^2+(w_2-1)^2 \quad \mathrm{s.t.}\ w_1 \geq 0,\  w_2 \geq 0,\ w_1w_2\leq 0.
	\end{align*}
	This example satisfies the MPCC-LICQ at all feasible points.	
	The point $w^*=(0,0)$ is C-stationary with the multipliers $\xi = -2$, $\nu = -2$. 
	However, it has two descent directions $d = (0,1)$ and $d = (1,0)$, and is thus not B-stationary.
	The points $w = (1,0)$ and $w =(0,1)$ are local minimizers and S-stationary points.
\end{example}
Arguably, stationarity conditions that allow first-order descent directions might be considered as too weak. 
S-stationarity is the only multiplier-based and non-combinatorial condition that has a chance to be equivalent to B-stationarity.
S-stationarity corresponds to the KKT conditions of the RNLP and it implies B-stationarity because of the inclusion relation in Eq. \eqref{eq:mpcc_feasible_sets}. 
Conversely, if the MPCC-LICQ holds, then B-stationarity implies also S-stationarity~\cite{Scheel2000}. 
Unfortunately, weaker conditions, e.g. the MPCC-MFCQ are already not sufficient for S-stationarity, see \cite[Example 3]{Scheel2000} for a counterexample.
In particular, M-stationarity is the strongest necessary condition under MPCC-MFCQ~\cite{Scheel2000,Flegel2006}.
\color{black}
This means that B-stationary points that are not S-stationary can not be identified via any stationarity concept based on the auxiliary NLPs~\cite{Veelken2009}. 
Instead, an exponential number of linear programs \eqref{eq:b_stationariry} must be solved for verification.

To summarize, under suitable CQs, local optimality is sufficient for any of the stationarity concepts defined here. 
However, everything weaker than S-stationarity is often considered to be too weak, since such points may allow first-order descent directions.
The problem in practice is that iterative MPCC methods are often attracted by M- or C-stationary points, even though the MPCC may have B- or S-stationary points~\cite{Leyffer2007}. 
There exist also second-order optimality conditions tailored to MPCCs. 
They are defined in terms of S-stationary points and the corresponding MPCC multipliers.
We omit their statement for brevity and refer the reader to \cite{Ralph2004,Scheel2000} for more details.

\section{MPCC solution methods}\label{sec:mpcc_methods}
In the last three decades, many tailored MPCC solution methods have been proposed. 
A recent survey of MPCC methods is given in \cite{Kim2020}. 
The references \cite{Fukushima2004,Hoheisel2013,Kanzow2015} provide comparisons of several solution strategies. 
We classify the MPCC solution methods into two distinct classes:
\begin{enumerate}
	[(a)]
	\item active-set-based/combinatorial methods,
	\item relaxation and smoothing-based methods.
\end{enumerate}
We review both classes, but give more details for the second class, because we will later benchmark them on the OCP-based problem set.

\subsection{Active-set based MPCC methods}
These methods explicitly treat the combinatorial nature of the complementarity constraints \eqref{eq:mpcc_cc}. 
They have the strongest convergence properties as they are usually guaranteed to converge to B-stationary points. 
They can be subdivided into branching methods \cite{Kim2020,Bard1990}, and active-set methods \cite{Leyffer2007,Kirches2022,Kim2020,Giallombardo2008,Guo2022}.
These methods rely on guessing the correct active set by solving a linear program with complementarity constraints (LPCC), such as \eqref{eq:b_stationariry} \cite{Leyffer2007,Kirches2022,Guo2022}.
If $d = 0$ solves the LPCC subproblem, a B-stationary point is found. 
To promote faster convergence rates, after fixing the active set, an equality-constrained QP can be solved~\cite{Leyffer2007,Kirches2022}. 
As the solution of an LPCC can be computationally expensive, Kirches et al. \cite{Kirches2022} regard LPCC only with complementarity and bound constraints. 
They suggest treating the remaining equality and inequality constraints in an augmented Lagrangian fashion. 
This is later done by Guo and Deng \cite{Guo2022}, where convergence to M-stationary points is proven.
The main practical drawback of this method class is the lack of robust open-source implementations. 
In contrast to the next group we treat, they are not easily implementable by using an available NLP code.
For more references we refer the reader to \cite[Section 5.6.2]{Lenders2018} and \cite[Section 3.2]{Kim2020}.

\subsection{Relaxation and smoothing-based methods}~\label{sec:relaxation_methods} 
The main idea behind relaxation-based methods is to replace \eqref{eq:mpcc_cc} with more regular constraints:
\begin{subequations}\label{eq:relaxation_general}
\begin{align}
	&G_i(w) \geq 0,\ H_i(w)\geq 0,\ &\text{for all}\ i \in \{1,\ldots, m\},  \\
	&\Phi(G_i(w),H_i(w),\sigma) \leq 0,\ &\text{for all}\ i \in \{1,\ldots, m\}. \label{eq:relaxation_general_deg}
\end{align}
\end{subequations}
which does not lead to the violation of LICQ and MFCQ. 
The function $\Phi: \R \times \R \times \R\to \R$ is a regularization function.
A smaller value of $\sigma$ yields a better approximation to the original problem and for $\sigma = 0$ we recover the original constraint \eqref{eq:mpcc_cc}.
Next one solves a (finite) sequence of these NLPs, by driving the parameter $\sigma> 0$ to zero. 
The availability of robust NLP codes makes their implementation easy and practical.
We denote the solution of the initial MPCC \eqref{eq:mpcc} by $w^*$ and the solution of the regularized NLP by $w^*(\sigma)$. 
The obvious goal is that $w^*(\sigma) \to w^*$ as $\sigma \to 0$.
Hoheisel et al. \cite{Hoheisel2013} provide a detailed numerical and theoretical comparison of several methods from this family.
If the problems are solved exactly, under mild assumptions, accumulation points of the sequence of solutions $w^*(\sigma_k)$ are usually C-stationary points~\cite{Scholtes2001,Ralph2004,Hoheisel2013}. 
Solving the NLPs inexactly usually weakens the convergence results \cite{Kanzow2015}. 
Of course, stronger assumptions also result in convergence to M- or S-stationary points. 
In the sequel, we discuss several examples of functions $\Phi$ and strategies for driving $\sigma$ to zero. 

\subsubsection{Direct solution}\label{sec:mpcc_direct}
One can use standard NLP techniques such as Sequential Quadratic Programming (SQP) and Interior Point (IP) methods for directly solving~\eqref{eq:mpcc}~\cite{Fletcher2002a,Fletcher2004,Leyffer2006a,Liu2004a}. 
Despite the degeneracy discussed so far, this approach can sometimes perform surprisingly well in practice. 
However, it tends to have convergence difficulties or to converge to spurious stationary points if the MPCC-LICQ does not hold.
Fletcher and Leyffer study the practical performance of SQP methods on numerous MPCCs in \cite{Fletcher2002} and investigate their local convergence properties in \cite{Fletcher2004}.
Under the assumptions that the MPCC-LICQ and that all QPs remain feasible, and other technical assumptions, they show quadratic convergence to S-stationary points. 
Interior-point methods perform reasonably well if applied directly to the NLP formulation \eqref{eq:mpcc}~\cite{Thierry2020}. 
However, their performance improves when paired with relaxation and exact penalty formulations, as we will highlight several times below.

\paragraph{NLP formulations with NCP functions}
In this approach the complementarity constraints \eqref{eq:mpcc_cc} are replaced by a C-function $\Phi_{\mathrm{C}}(G(w),H(w)) = 0$. 
The resulting NLP is solved by a standard globalized NLP solver. 
The use of SQP methods in such formulations was studied by Leyffer~\cite{Leyffer2006a}. 
If the used C-function is not differentiable at $(0,0)$, its subgradient can be used~\cite{Leyffer2006a}.
Evidently, using squared versions (or higher powers) of C-functions will not improve the situation and lead to a violation of LICQ at $(0,0)$ since we obtain a zero-gradient at this point.
Similar to \cite{Fletcher2004}, assuming MPCC-LICQ, Lipschitz continuity of the NLP functions and their derivatives, and other technical assumptions, local superlinear convergence to S-stationary points was shown. 
Leyffer~\cite{Leyffer2006a} tests this approach on a wide number of test problems and shows that different NCP functions can lead to large differences in performance. 

\subsubsection{Global relaxation/smoothing method by Scholtes}\label{sec:scholtes}
\begin{figure}[t]
	\centering
	{\includegraphics[width=\textwidth]{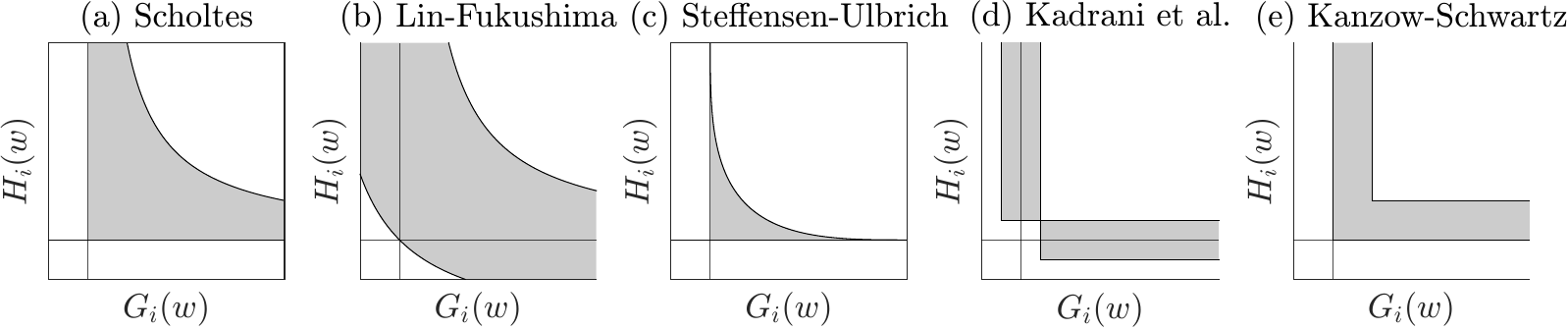}}
	\caption{Illustration of the regularized complementarity sets.}
	\label{fig:mpcc_relaxation}
\end{figure}

Scholtes introduced  arguably the easiest-to-implement approach and relaxes the complementarity constraint by using~\cite{Scholtes2001}:
\begin{align}\label{eq:scholtes}
	\Phi_{\mathrm{S}}(G_i(w),H_i(w),\sigma) = G(w)_i H_i(w) - \sigma.
\end{align}
An illustration of the relaxed feasible set is given in Figure \ref{fig:mpcc_relaxation} (a).
Alternatively, the bilinear term might be smoothed by using $\Phi_{\mathrm{S}}(G_i(w),H_i(w),\sigma) = 0$ in \eqref{eq:relaxation_general_deg}.
Lumped versions $\Phi_{\mathrm{S}}(G(w),H(w),\sigma)  = G(w)^\top H(w) - \sigma$ are also used frequently~\cite{Scholtes2001}. 
Observe that in contrast to the smoothed variant, the relaxed version contains the feasible set $\Omega_{\mathrm{MPCC}}$, and one might find with it a stationary point without driving $\sigma \to 0$.
Assuming MPCC-LICQ, Scholtes \cite{Scholtes2001} shows convergence to C-stationary points. 
Hoheisel et al. \cite{Hoheisel2013} obtain the same result under the weaker MPCC-MFCQ.
Ralph et al. \cite{Ralph2004} study the convergence speed of this approach and show that, under rather strict assumptions, the local solution map $w^*(\sigma)$ of the relaxed formulation is a piecewise continuous function and that the solution converges with a rate $O(\sigma)$. 
Milder assumptions result in the rate $O(\sigma^{\frac{1}{2}})$ for the relaxed variant and $O(\sigma^{\frac{1}{4}})$ for the smoothed variant.

For more efficacy, Raghunathan and Biegler \cite{Raghunathan2005}, Liu and Sun \cite{Liu2004a} propose interior-point methods where the relaxation parameter $\sigma$ is proportional to the barrier parameter $\tau$. 
Under some stronger assumptions, this strategy is shown to converge to S-stationary points.

\paragraph{Smoothed NCP functions}
Early MPCC algorithms considered smoothed and everywhere differentiable variants of C-functions. 
These methods are closely related to Scholtes' approach since one can often by simple algebraic manipulation obtain the same feasible set as with using \eqref{eq:scholtes}.
In our experiments, we will test three different smoothed NCP functions:
the smoothed Fischer-Burmeister (FB) function ${\Phi}_{\mathrm{FB}}(a,b,\sigma) = a+b-\sqrt{a^2+b^2+\sigma^2}$, 
the Natural Residual (NR) function ${\Phi}_{\mathrm{NR}}(a,b,\sigma) = (a+b-\sqrt{(a-b)^2+\sigma^2})$, 
the Chen-Chen-Kanzow (CCK)~\cite{Chen2000a} function ${\Phi}_{\mathrm{CCK}}(a,b,\sigma) = \lambda(a+b-\sqrt{a^2+b^2+\sigma^2})+(1-\lambda)(ab-\sigma)$, where $\lambda\in (0,1)$ is a parameter, which is set to 0.5 in our implementations.
Facchinei et al. \cite{Facchinei1999} show the convergence of this approach to C-stationary points.
In our implementation, we use the smoothed versions of the NCP functions in the form of \eqref{eq:relaxation_general}. 
Hence, we obtain relaxations of the original problem~\eqref{eq:mpcc}.

\subsubsection{The smooth relaxation method by Lin and Fukushima} 
This method is similar to Scholtes' regularization and replaces the complementarity conditions \eqref{eq:mpcc} by:
\begin{align*}
	&G_i(w) H_i(w) \leq \sigma^2, &\ \text{for all}\ i \in \{1,\ldots,m\} ,\\
	&(G_i(w)+\sigma)( H_i(w) + \sigma) \geq \sigma^2, &\ \text{for all}\ i \in \{1,\ldots,m\}.
\end{align*}
Figure \ref{fig:mpcc_relaxation} (b) shows an illustration of the feasible set. 
Lin and Fukushima \cite{Lin2005} obtain similar convergence results as Scholtes \cite{Scholtes2001}. 
Hoheisel et al. \cite{Hoheisel2013} extend this result by proving convergence to C-stationary points under the MPCC-MFCQ. 
Moreover, they show that the feasible points of the relaxed NLP satisfy the MFCQ in a neighborhood of a point $x \in \Omega_{\mathrm{MPCC}}$.

We continue with reviewing more sophisticated regularization schemes that converge to M-stationary points under fairly mild conditions~\cite{Hoheisel2013}. 
As we will see later, this may not always imply better performance in practice.
\subsubsection{The local relaxation method by Steffensen and Ulbrich}\label{sec:steffensen_ulbrich}
Almost all regularization methods make global changes to the feasible set. 
Motivated by the fact that most difficulties arise for degenerate complementarity pairs $i \in \mathcal{I}_{00}$, Steffensen and Ulbrich follow a different approach \cite{Steffensen2010,Veelken2009}. 
Their main idea is to relax the complementarity constraint only locally at the corner of the L-shaped set arising from the complementarity constraints. 

The relaxation is achieved by the following steps: the L-shaped set is rotated with a linear transformation by $\frac{\pi}{4}$ counterclockwise for every complementarity pair, and one obtains the graph of the abs-function. 
On the interval $[-\sigma,\sigma]$, the kink is replaced by a sufficiently smooth function such that the continuity of the functions and their derivatives is preserved at the interval boundaries. 
Finally, the inverse transformation is carried out, and a locally relaxed set is obtained, cf. Figure \ref{fig:mpcc_relaxation} (c). 
This reasoning expressed in equations reads as
\begin{align*}
	&\Phi_{\mathrm{SU}}(G_i(w),H_i(w),\sigma) \leq 0, & \ \text{for all}\ i \in \{1,\ldots,m\},
\end{align*}
where $\Phi_{\mathrm{SU}}: \R \times \R \times \R \to \R$ is defined in terms of the auxiliary functions  $\phi^{a}_{\mathrm{SU}}: \R \times \R \to \R$ and $\phi^{b}_{\mathrm{SU}}: [-1,1] \to \R$ as follows
\begin{align*}
	\Phi_{\mathrm{SU}}(y_1,y_2;\sigma) = y_1 + y_2 - \phi_{\mathrm{SU}}^{a}(y_1-y_2,\sigma), \\
	\phi_{\mathrm{SU}}^{a}(z,\sigma)  = \begin{cases}
		|z|, &\textrm{ if } |z| \geq \sigma,\\
		\sigma \phi_{\mathrm{SU}}^{b}(\frac{z}{\sigma}),  &\textrm{ if } |z| < \sigma.
	\end{cases}
\end{align*}
The function $\phi_{\mathrm{SU}}^{b}$ has to satisfy some smoothness and mononoticitiy properties~\cite{Steffensen2010}.
For our experiments we implement two variants of such functions proposed in \cite{Steffensen2010}: $\phi_{\mathrm{SU}}^b(z) = \frac{1}{8} (-z^4+6z^2+3)$ and $\phi_{\mathrm{SU}}^b(z) = \frac{2}{\pi} \sin(z \frac{\pi}{2} +\frac{3\pi}{2})+1$ .
\color{black}
Under the MPCC-CRCQ (cf. \cite[Defintion 2.4]{Hoheisel2013}) convergence to C-stationary and under the MPCC-LICQ to M-stationary points is shown \cite{Steffensen2010}.
\subsubsection{The nonsmooth relaxation method by Kadrani et al.}
Another interesting relaxation by Kadrani et al. \cite{Kadrani2009} reads as:
\begin{align*}
	&G_i(w)  \geq -\sigma,\; H_i(w)  \geq -\sigma, & \ \text{for all}\ i \in \{1,\ldots,m\},\\
	&(G_i(w)-\sigma)( H_i(w) - \sigma) \leq 0,  &\ \text{for all}\ i \in \{1,\ldots,m\}.
\end{align*}
Figure \ref{fig:mpcc_relaxation} (d) illustrates the nonsmooth feasible set obtained from the constraint above. 
The convergence study of \cite{Kadrani2009} is carried out assuming the MPCC-LICQ. 
Once again, Hoheisel et al. \cite{Hoheisel2013} improve the result and show convergence to M-stationary points under the MPCC-CPLD (cf. \cite[Defintion 2.4]{Hoheisel2013}, a CQ weaker than the MPCC-MFCQ and stronger than the MPCC-ACQ). 
It is evident from the structure of the feasible set of this relaxation that verifying standard CQs is more difficult. 

\subsubsection{The relaxation method by Kanzow and Schwartz}
This relaxation has stronger theoretical properties than the previous one and a more satisfactory shape of the feasible set, cf. Figure \ref{fig:mpcc_relaxation} (e)~\cite{Kanzow2013}. 
In contrast to the approach of Kadrani et al., it contains the feasible set of the MPCC. 
This relaxation is modeled with the following equations:
\begin{align*}
	&\Phi_{\mathrm{KS}}(G_i(w),H_i(w),\sigma) \leq 0,  &\textit{for all}\ i \in \{  1,\ldots,m\},
\end{align*}
with $\Phi_{\mathrm{KS}}: \R \times \R \times \R \to \R$ and $\phi_{\mathrm{KS}}: \R \times \R \to \R$ where  $\Phi_{\mathrm{KS}}(y_1,y_2,\sigma) = \phi_{\mathrm{KS}}(y_1-\sigma,y_2-\sigma)$ and 
\begin{align*}
	\phi_{\mathrm{KS}}(y_1,y_2)  &= \begin{cases}
		y_1 y_2, &\textrm{ if } y_1+y_2 \geq 0,\\
		-\frac{1}{2}(y_1^2+y_2^2), &\textrm{ if } y_1+y_2 <0.	
	\end{cases}
\end{align*}
The function $\phi_{\mathrm{KS}}$ is a continuously differentiable C-function \cite{Kanzow2013}.
Under the MPCC-CPLD (cf. \cite[Definition 2.4]{Hoheisel2013}) convergence to M-stationary points is shown \cite{Kanzow2013,Hoheisel2013}. 

\subsubsection{Exact penalty methods}\label{sec:mpcc_penalty}
Exact penalty reformulations are one of the most often used approaches to treat degenerate NLPs~\cite{Benson2006,Byrd2012,Gill2005}.
In exact penalty algorithms for MPCCs, the bilinear term \eqref{eq:mpcc_cc} is added to the objective in some suitable form and multiplied by a penalty factor $\rho$. 
To be consistent with our notation above, and the implementations in \nosnoc~\cite{Nurkanovic2022b}, we use $\rho^k = \frac{1}{\sigma_k}$. 
Assuming sufficient regularity of other constraints and having the bilinear term in the objective, we obtain a regular NLP.
Under suitable assumptions and for a sufficiently large and finite $\rho$, the solution $w^*(\sigma)$ matches the solution $w^*$ of the initial MPCC after a single NLP solve~\cite{Ralph2004}. 
In practice, a sequence of NLPs is solved to improve the convergence and to estimate the correct penalty parameter value. 
One of the simplest formulations is to penalize the term $G(w)^\top H(w)$ in the objective. 
This corresponds to the $\ell_1$ norm of the complementarity residual. 
Therefore, we solve a sequence of the following NLPs:
\begin{subequations}\label{eq:mpcc_ell_1_penalty}
	\begin{align}
		\min_{x} \quad &   f(w) + \frac{1}{\sigma_k}  G(w)^\top H(w)\\
		\mathrm{s.t.} \quad  &g(w) = 0 ,\\
		&h(w)\geq 0, \\ 
		&G(w) \geq 0,\ H(w) \geq 0.
	\end{align}
\end{subequations}
This approach was first proposed in \cite{Ferris1999} for solving practical problems. 
Anitescu \cite{Anitescu2000a} provided the first convergence analysis for the $\ell_1$ penalty approach paired with active-set SQP methods. 
Ralph et al. \cite{Ralph2004} show that an MPCC solution $w^*$ is also a solution to \eqref{eq:mpcc_ell_1_penalty} for a sufficiently large $\rho$ and that regularity conditions of the MPCC (e.g., MPCC-LICQ) imply regularity of the NLP \eqref{eq:mpcc_ell_1_penalty}.
However, if the local minimizers are only B-stationary but not S-stationary points, the penalty parameter must grow to infinity \cite{Kim2020}.

Leyffer et al. \cite{Leyffer2006} propose an interior-point algorithm to solve the NLP \eqref{eq:mpcc_ell_1_penalty} while dynamically updating the penalty parameter $\rho$. 
For each fixed $\rho^k$, the barrier subproblem is solved inexactly to a tolerance proportional to the barrier parameter $\tau^k$. 
Strategies to steer the penalty parameter that avoid too large increases and unbounded subproblems are proposed as well. 
Fukushima et al. \cite{Fukushima1998} suggest an SQP method paired with a penalized NCP function.
Hu and Ralph \cite{Hu2004} relate relaxation methods of \cite{Scholtes2001} and give conditions for convergence to B-stationary points. 
They study more general formulations than \eqref{eq:mpcc_ell_1_penalty} and suggest for example to use $\sum_{i=1}^{m} \Phi_{\mathrm{FB}}(G_i(w),H_i(w))^3$ as a penalty function.
Furthermore, by comparing the KKT conditions, Leyffer et al. \cite{Leyffer2006} show that there exists a one-to-one correspondence between the iterates of the penalty and a smoothing Scholtes approach (i.e., the bilinear constraint  in \eqref{eq:mpcc_cc} is replaced by $G_i(w)H_i(w) = \sigma_k$ ).

Hall et al. \cite{Hall2022,Hall2023} propose a Sequential Convex Programming (SCP) method for solving the penalty problem arising from quadratic programs with complementarity constraints. 
In particular, the method makes use of the fact that QP matrices do not need to be re-factorized in an SCP approach, which enables fast and cheap iterations. 
The algorithm is paired with an exact analytic line search.

The great practical difficulty in exact penalty methods is steering the penalty parameter. 
Byrd et al. \cite{Byrd2012} introduce a line-search SQP  $\ell_1$-exact penalty method for general degenerate NLPs. 
They propose penalty update rules based on solving LPs and QPs with a trust region to predict the decrease of the merit function. 
Favorable theoretical properties and good numerical performance on a series of test problems including MPCCs are reported. 
Thierry and Biegler \cite{Thierry2020} adapt the $\ell_1$ strategy of Byrd et al. \cite{Byrd2012} including the penalty steering rules, to solve degenerate problems with \ipopt~\cite{Waechter2006}. 
Good practical performance is reported on the \MacMPEC\ test set \cite{Leyffer2000} with an improvement in terms of speed and robustness compared to a direct application of \ipopt.

Moreover,  using an $\ell_\infty$  penalty function for MPCC is also very common \cite{Benson2006}. 
The complementarity constraints \eqref{eq:mpcc_cc} are replaced by:
\begin{align*}
	&G(w) \geq 0,\ H(w) \geq 0,\ G_i(w) H_i(w) \leq s, \;\text{for all}\ i\in \{1,\ldots, m\},
\end{align*}
where $s\in \R$ is a slack variable, which may have an upper bound $\bar{s}>0$.
The term $\frac{1}{\sigma} s$ is added to the objective function.
This enables us to express the $\ell_\infty$ norm of the complementarity constraints residual smoothly. 
{Similarly}, if the bilinear terms are lumped together, and we use the constraint  $G(w)^\top H(w) \leq s$, we end up with an $\ell_1$ formulation that is equivalent to \eqref{eq:mpcc_ell_1_penalty}.

A mixture of the approaches above is the elastic mode, which takes an $\ell_1$ norm of the bilinear complementarity terms and an $\ell_\infty$ penalization of the relaxation of the standard equality and inequality constraints. 
Anitescu et al. \cite{Anitescu2007} show under the MPCC-LICQ and other assumptions the global convergence of an elastic mode SQP approach to C-, M- and S-stationary points. 
The elastic mode with a fixed penalty parameter is implemented in SNOPT as a fallback strategy if an infeasible or unbounded QP is detected~\cite{Gill2015}. 

Finally, we mention the family of augmented Lagrangian methods for MPCC, which also belong to the class of penalty methods \cite{Fukushima2004,Izmailov2012a}. 
We do not treat them in detail here. 
Assuming MPCC-LICQ, and that the sequence of Lagrange multipliers is bounded, convergence to S-stationary points is shown in \cite{Fukushima2004,Izmailov2012a}. 

\subsubsection{Lifting methods}\label{sec:mpcc_lifting}
Lifting methods are somewhat in between relaxation and penalty methods. 
The main idea is to introduce \textit{lifting variables} and regard a more regular feasible set in a higher-dimensional space whose orthogonal projection is the L-shaped set, coming from the complementarity constraints.
Some of them require penalization of the lifting variables to recover the solution of the initial problem \cite{Hatz2013}, and others might require additional regularization \cite{Stein2012a}. 
Unfortunately, they have weaker theoretical properties than regularization methods, cf. Section \ref{sec:mpcc_summary}.
Thus, we do not implement these methods and do not treat them in further detail.
We mention the methods of Stein \cite{Stein2012a}, Hatz et. al \cite{Hatz2013}, and Izmailov et al.\cite{Izmailov2012a}.

\subsection{Steering the homotopy parameter to zero}\label{sec:steering}
The initial value for the homotopy parameter $\sigma_0$ and deciding how to steer it to zero play an important role in the practical performance of the relaxation-based methods.
In our implementations, we take three different approaches in steering the relaxation in \eqref{eq:relaxation_general}:
\begin{enumerate}
	\item directly change the relaxation parameters $\sigma$ outside of the homotopy loop, which is the standard approach in most of the literature,
	\item steer a single relaxation parameter with an $\ell_\infty$ penalty approach to zero,
	\item steer several relaxation parameters with an $\ell_1$ penalty approach to zero.
\end{enumerate}
The ways in which the different approaches manifest in a relaxed NLP are summarized in Table~\ref{tab:parameter_steering}.
In the standard approach, we simply use a fixed parameter $\sigma_k$ for \eqref{eq:relaxation_general_deg} in every NLP solve, and update the parameter after every iteration via:
\begin{align}\label{eq:parameter_homotopy}
	\sigma_{k+1} = \kappa \sigma_k,
\end{align}
with $\kappa \in (0,1)$. 
Alternatively, we may use the update formula (as often used in IP methods \cite{Nocedal2006}):
\begin{align*}
	\sigma_{k+1} = \min(\kappa \sigma_k,\sigma_k^\eta),
\end{align*}
with $\eta >1$.
Next, we may let the optimizer steer the relaxation, by using e.g., the same scalar variable $s$ in all \eqref{eq:relaxation_general_deg}. 
The slack variable is pushed to zero by being more and more penalized in the objective function, with a weight of $\rho_k = \frac{1}{\sigma_k}$.
In the case of Scholtes' relaxation, we recover the $\ell_\infty$ penalty approach discussed in Section \ref{sec:mpcc_penalty}.
Similarly, we can add a new scalar variable $s_i$ for every constraint in \eqref{eq:relaxation_general_deg} and penalize their deviation from zero in the objective and weighted by $\rho_k = \frac{1}{\sigma_k}$.
In the case of Scholtes' relaxation, we recover the exact $\ell_1$ penalty approach and obtain a NLP equivalent to~\eqref{eq:mpcc_ell_1_penalty}.

\begin{table}[b]
	\centering
	{\def\arraystretch{1.5}\tabcolsep=5pt
		\begin{tabular}{|c|c|c|}
			\hline 
			{Standard relaxation}  & $\ell_{\infty}$ penalty relaxation & $\ell_1$ penalty relaxation\\
			\hline
			$
			\begin{aligned}
				\underset{w\in \R^{n_w}}{\mathrm{min}} \;  & f(w) \\
				\textrm{s.t.} \quad 
				0 &= g(w),\\
				0 &\leq h(w),\\
				0 & \leq G(w),\\
				0 & \leq H(w),\\
				0&\geq\Phi (G_i(w),H_i(w),\sigma_k),\\
				& \ \text{for all}\ i \in \{1,\ldots,m\}.
			\end{aligned}
			$
			&
			$
			\begin{aligned}
				\underset{\substack{w\in \R^{n_w}\\ s\in \R}}{\mathrm{min}} \;  & f(w) + \frac{1}{\sigma_k} s  \\
				\textrm{s.t.} \quad 
				0 &= g(w)\\
				0 &\leq h(w)\\
				0 & \leq G(w)\\
				0 & \leq H(w)\\
				0&\geq\Phi (G_i(w),H_i(w),s),\\
				& \ \text{for all}\ i \in \{1,\ldots,m\}.
			\end{aligned}
			$
			&
			$
			\begin{aligned}
				\underset{\substack{w\in \R^{n_w}\\ s\in \R^{m}}}{\mathrm{min}} \;  & f(w) + \frac{1}{\sigma_k} \sum_{i=1}^{m} s_i  \\
				\textrm{s.t.} \quad 
				&0 = g(w),\\
				&0 \leq h(w),\\
				&0  \leq G(w),\\
				&0  \leq H(w),\\
				&0\geq\Phi (G_i(w),H_i(w),s_i),\\
				& \ \text{for all}\ i \in \{1,\ldots,m\}.
			\end{aligned}
			$
			
			\\
			\hline
		\end{tabular}
		\normalsize}
	\caption{Homotopy parameter steering strategies.}
	\label{tab:parameter_steering}
\end{table}
In our experiment, we solve all NLPs in the homotopy loop to the same prescribed tolerance. 
Alternatively, one may solve the NLPs in the homotopy loop inexactly.
Some authors suggest updating the relaxation parameter simultaneously or tying it to other relaxation parameters of the algorithm.
Raghunathan and Biegler \cite{Raghunathan2005} update the parameter $\sigma$ simultaneously with the barrier parameter $\tau$ in an IP method.
Similarly, in Leyffer et al. \cite{Leyffer2006} in an IP-$\ell_1$ exact penalty approach the parameter update is related to the update of $\tau$.
Lin and Othsuka \cite{Lin2022} use a non-interior point method with Scholtes' relaxation. There the parameter of the Scholtes' relaxation and the complementarity constraint relaxation parameter in the KKT conditions are updated simultaneously.

To make the various relaxations better comparable and less dependent on the scaling of $\sigma$, we bring them to a similar scale for a given $\sigma$.
In particular, for a given $\sigma$, we want that the distance of $\Phi(a,b,\sigma)  =0$ to the origin is always the same, independent of the particular choice of $\Phi$.
Moreover, for a given $\kappa$ we require that this distance shrinks with the same rate for all relaxations.
By simple algebraic manipulations one can find how $\sigma$ needs to enter a given relaxation function. 
Our initial experiments confirm that this makes the performance more consistent and better comparable.

\subsection{Summary of MPCC methods}\label{sec:mpcc_summary}
	\begin{table}[t]
		\scriptsize
		{\def\arraystretch{1.5}\tabcolsep=5pt
			\begin{tabular}{|c|c | p{2.24cm} p{1.0cm} p{1.9cm} p{1.1cm}|}
				\hline 
				\textbf{{Type}} &\textbf{{Method}}  &{CQ for MPCC} & {Limiting Point}  &{Subproblem NLP satisfies} & {Citation} \\
				\hline
				Direct& Fletcher et al.	    & MPCC-LICQ &  S & GCQ & \cite{Fletcher2006,Fletcher2004}\\
				& Leyffer    			& MPCC-LICQ &  S  & GCQ & \cite{Leyffer2006}\\
				\hline 
				Regular- &Scholtes 	        & MPCC-MFCQ &  C  & MFCQ &\cite{Scholtes2001,Hoheisel2013}\\
				ization &Lin-Fukushima       & MPCC-MFCQ &  C  & MCFQ &\cite{Lin2005,Hoheisel2013}\\
				&Kadrani et al. 		& MPCC-MFCQ &  M  & GCQ &\cite{Hoheisel2013,Kadrani2009}\\
				&Steffensen-Ulbrich  & MPCC-CPLD &  C  & ACQ &\cite{Hoheisel2013,Steffensen2010}\\
				&Kanzow-Schwartz     & MPCC-CPLD &  M  & GCQ &\cite{Hoheisel2013,Kanzow2013}\\
				&Raghunathan-Biegler & MPCC-LICQ &  {C}  &  MFCQ &\cite{Raghunathan2005}\\
				\hline
				Lifting&Stein    	        & MPCC-LICQ &  C  & LICQ &\cite{Stein2012a}\\
				&Izmailov-Solodov    & MPCC-LICQ &  C & LICQ &\cite{Izmailov2012b}\\
				&Hatz et al.  	    & MPCC-LICQ &  {S}  & GCQ &\cite{Hatz2013}\\
				\hline
				Penalty&$\ell_1$-Penalty    & MPCC-LICQ &  {S}  & LICQ &\cite{Ralph2004,Leyffer2006a}\\
				&Leyffer et al.   & MPCC-LICQ &  C   & LICQ &\cite{Leyffer2006a,Ralph2004}\\
				&$\ell_\infty$-Penalty  & MPCC-LICQ &  {S}   & LICQ &\cite{Ralph2004,Benson2006}\\
				&Elastic mode	    & MPCC-LICQ &  C  & MFCQ &\cite{Anitescu2005,Anitescu2007}\\
				\hline
				Active-set &Leyffer-Munson    & MPCC-MFCQ  &  B  & GCQ &\cite{Leyffer2007}\\
				&Kirches et al.    & -  &  B  & GCQ &\cite{Kirches2022}\\
				&Guo-Deng    & MPCC-RCPLD  &  M  & GCQ &\cite{Guo2022}\\
				\hline
			\end{tabular}
			\normalsize}
		\caption{Overview of convergence properties of MPCC methods.}
		\label{tab:mpcc}
		\vspace{-0.5cm}
	\end{table}
	
	With a standard NLP solver at hand, all methods from Section \ref{sec:relaxation_methods} paired with the steering strategies from Section \ref{sec:steering} are easy to implement.
	Table \ref{tab:mpcc} provides an overview of known convergence results for direct, regularization, lifting, penalty methods, and active-set methods.
	Together with the diagram in Figure \ref{fig:mpcc_relations}, it can help one to decide which algorithm to choose to compute a stationary point of the MPCC~\eqref{eq:mpcc}. 
	The strongest multiplier-based stationarity concept is S-stationarity, followed by M-, C-, A- and W-stationarity, sorted from stronger to weaker. 
	One should not be discouraged by the weaker limiting points of the relaxation methods. 
	In contrast to the other methods, these results are obtained under much weaker assumptions. 
	In general, more restrictive assumptions give a better result. 
	For example, the global relaxation method by Scholtes \cite{Scholtes2001} converges to B-stationary points under the MPCC-LICQ and the \text{upper-level strict complementarity}, cf. \cite [Definition 2.6]{Ralph2004}.
	We note that all methods converging to an S-stationary point under the MPCC-LICQ also assume other restrictive assumptions, such as the upper-level strict complementarity \cite{Fletcher2006,Hatz2013,Leyffer2006a,Raghunathan2005}.
	
	In practice, the NLP subproblems cannot be solved exactly, due to the solver tolerances and finite arithmetic precision.
	However, solving the NLPs in the sequence inexactly can weaken the convergence results~\cite{Kanzow2015,Anitescu2007}. 
	For example, under inexact solves the methods of Kadrani et al., Kanzow-Schwartz \cite{Kanzow2013} and Steffensen-Ulbrich \cite{Steffensen2010} converge only to W-stationary points.  
	Surprisingly, the methods of Scholtes and Lin-Fukushima are immune to this, and they still converge to C-stationary points~\cite{Kanzow2015}.
	However, if some stronger assumptions are not satisfied, usually the MPCC-LICQ, they can experience slow convergence rates and converge to weaker points. 
	This motivated the development of combinatorial methods, which have excellent theoretical properties but currently no mature open-source implementations.

\section{A benchmark set of MPCC from nonsmooth OCPs}\label{sec:nosbench}
The main contribution of this paper is the introduction of a benchmark suite of MPCCs that come from nonsmooth optimal control and simulation problems.
We discuss the problem format, provide references for the original continuous time OCP and simulation problems, and discuss how we generate MPCCs from them.
Moreover, we split the whole problem set into several subsets to facilitate the testing of the variety of available algorithms.
\subsection{Problem format}
In the benchmark we provide the MPCCs in a slightly different format than in Eq. \eqref{eq:mpcc}:
\begin{subequations}\label{eq:mpcc_nosbench}
	\begin{align}
		\underset{w\in \R^{n}}{\mathrm{min}} \; &f(w,p)\\
		\textrm{s.t.} \quad 
		&\ell_w \leq w \leq u_w,\\
		&\ell_g \leq  g (w,p) \leq u_g,\\
		&0 \leq G(w,p)  \perp H(w,p) \geq 0 \label{eq:mpcc_nosbench_cc},
	\end{align}
\end{subequations}
where $p \in \R^{n_p}$ a given parameter. 
If we disregard \eqref{eq:mpcc_nosbench_cc}, this format compiles with the \casadi\ NLP solver interface~\cite{Andersson2019}. 
All problem functions in \eqref{eq:mpcc_nosbench} are \casadi\ functions generated via \nosnoc. 
By treating the complementarity constraint with some of the methods described in Section \ref{sec:relaxation_methods} we obtain an NLP that can be directly passed to the \casadi\  NLP solver interface.
On the benchmarks homepage \footnote{\nosbench\ homepage: \nosbenchurl.} we provide all MPCCs in the form of Eq. \eqref{eq:mpcc_nosbench} format using a JSON object with the following components:
\begin{itemize}[itemsep=0.5em]
	\item $w$ as a string encoded \casadi\ variable with the key ``\texttt{w}'', along with the $\ell_{w}$ and $u_{w}$ with the labels ``\texttt{lbw}'' and ``\texttt{ubw}'', respectively,
	\item $f(w)$ as a string encoded \casadi\ function with the key ``\texttt{f\_fun}'',
	\item $g(w)$ as a string encoded \casadi\ function with the key ``\texttt{g\_fun}'', along with the $\ell_g$ and $u_{g}$ with the labels ``\texttt{lbg}'' and ``\texttt{ubg}'', respectively,
	\item $G(w)$ and $H(w)$ as string encoded \casadi\ functions with keys ``\texttt{G\_fun}'' and ``\texttt{H\_fun}'' respectively,
	\item Parameters $p$ and their values as string encoded \casadi\ variables and an array of doubles with keys ``\texttt{p}'' and ``\texttt{p\_val}'' respectively.
\end{itemize}
From this, a user can simply reconstruct the problem by loading in all of the components and using the provided \casadi\ deserialization functionality to interface their solver with \nosbench\  problems.

\subsection{Description of the benchmark collection}
The \nosbench\ problems are generated by using the algorithmic toolchain available in \nosnoc. 
In particular, we regard both simulation and optimal control problems, reformulate them, and discretize them with the FESD method \cite{Nurkanovic2022,Nurkanovic2023b}.
We list in Table~\ref{tab:NOSBENCH-ODEs} the origins of each of these problems in continuous time along with references to them in the literature. 
We briefly describe each of the systems along with a classification. 
Moreover, we also list the number of MPCCs we generated from the particular continuous-time problem. 
For more details on each discrete-time problem, we refer the reader to the benchmark's homepage.

\subsubsection{Reformulation and discretization options}
Depending on the reformulation and discretization method, the continuous-time problems can be reformulated into significantly different MPCCs. 
In order to generate problems of varying complexity and internal structure we vary several discretization and MPCC parameters. 
In all cases we obtain an MPCC of the form of Eq. \eqref{eq:ocp_discrete_time}. 
We list some variations, which are available in \nosnoc:
\begin{itemize}
	\item We distinguish between simulation and optimal control problems (OCPs).
	\item If the system has a nonsmooth {right hand side} (r.h.s.), but no state jumps we treat it as Filippov system, which one can reformulate either via Stewart's  \cite{Nurkanovic2022} or via the Heaviside step reformulation \cite{Nurkanovic2023c} into a DCS. 
	Afterwards, they are discretized with the corresponding FESD method. 
	\item Complementarity Lagrangian Systems (CLS) with state jumps can be either treated directly with FESD-J \cite{Nurkanovic2023b} or reformulated via time-freezing into a Filippov system~\cite{Nurkanovic2021,Nurkanovic2023a}. 
	Hybrid systems with hysteresis are always reformulated via time-freezing into Filippov systems.
\end{itemize}
Moreover, we can vary the discretization parameters:
\begin{itemize}[itemsep=0.5em]
	\item $N$ - the number of control intervals in the OCP. 
	The controls are piecewise constant after the discretization.
	This is always set to $N=1$ in the case of simulation problems.
	\item $N_{\mathrm{FE}}$: number of integration steps (finite elements) within a control interval.
	\item The number of stage points used by the selected Runge-Kutta (RK) method within each finite element in the FESD discretization is given by $n_{\mathrm{s}}$. 
	\item Choice of underlying RK method (defined by its Butcher tableau) in FESD discretization. 
	In our experiments, we regard the Radau-IIA or Gauss-Legendre schemes. 
	\item There are several choices for grouping the cross-complementarity conditions in the FESD discretization, cf. \cite{Nurkanovic2022}. 
	They provide different sparsity to the complementarity constraints that enforce switch detection.
\end{itemize}
Furthermore, in some cases, we vary the problem parameters. 
This is done to provide some variety in the benchmark problems and to mitigate the effects of pre-tuning which has been done on some of these examples in the corresponding references. 
It also allows for simulation problems to be tested both in cases where there are no switches and in cases where switches must be detected.

\begin{figure}[t]
	\centering
	\subfloat[Number of inequality constraints versus number of variables.]{\includegraphics[width=0.48\textwidth]{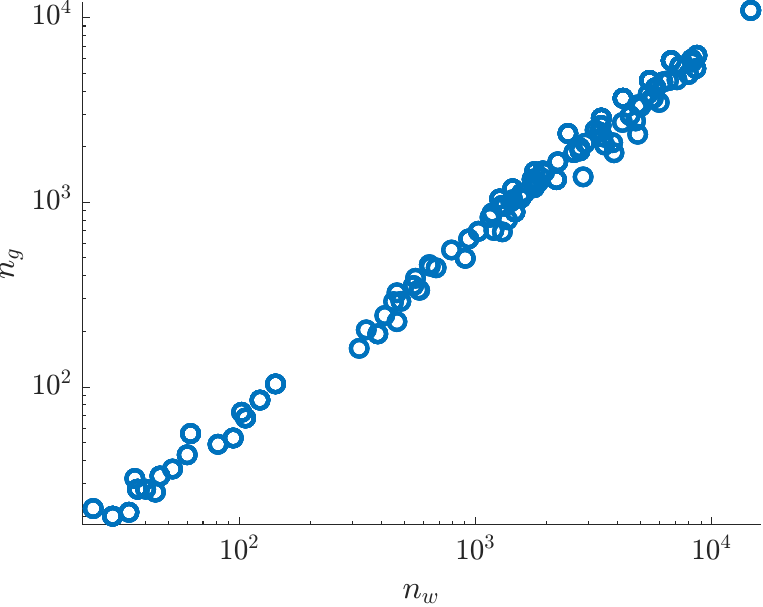}}
	\hfill
	\subfloat[Number of complementarity constraints versus number of variables.]{\includegraphics[width=0.48\textwidth]{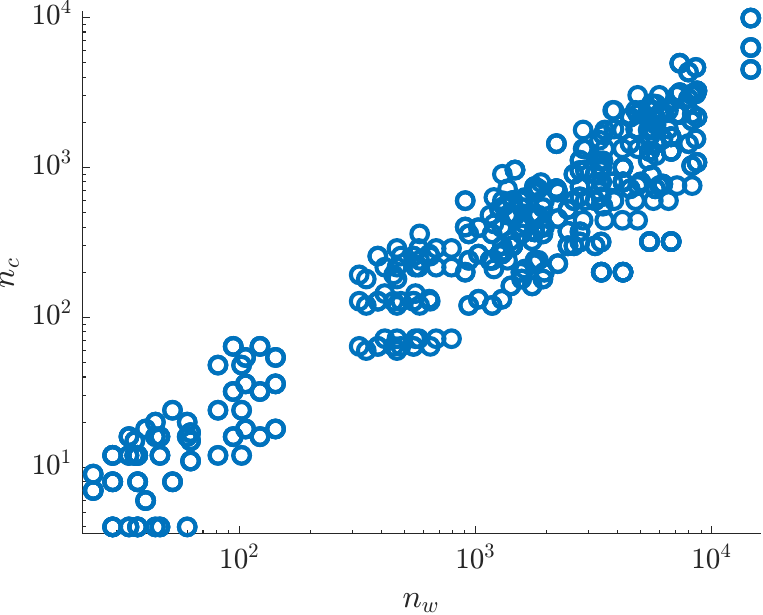}}
	\caption{Size characteristics of the \nosbench\ test set.}
	\label{fig:gvw-cvw}
\end{figure}

By varying all the aforementioned parameters we generate a total of {\NF} distinct MPCCs within the full problem set. 
Due to space limitations, we cannot provide full details in this paper, but they are available on the benchmark's homepage.
Figure \ref{fig:gvw-cvw} shows the number of primal variables versus the dimensions $n_g$ of $g(w)$ in \eqref{eq:mpcc_nosbench} and versus the number of complementarity pairs in each problem.

\subsubsection{Problem naming convention}
The names of the problem files are encoded with some information about the problem. 
This is done so that the formation of subsets of the whole benchmark can be done systematically.
The structure of these names is an underscore delimited string containing the following data in order: Problem name, parameter index, $N$, $N_{\mathrm{FE}}$,$n_{\mathrm{s}}$, RK method, DCS type (``Step'', ``Stewart'', or ``CLS''), cross complementarity mode, source, and a boolean flag if the problem is lifted into the vertical form.
The type of the problem is split into ``FIL'', Fillipov systems, ``IEC'', problems with only inelastic collisions, ``ELC'', problems with elastic (and possibly inelastic) collisions, and ``HYS'' for problems containing hysteresis.
For example the problem titled: \texttt{986EQ\_001\_001\_003\_2\_GL\_STEP\_7\_FIL\_1} would be the earthquake example from \cite{Calvo2016}, with the first parameter set and $N = 1$, $\NFE=3$, $n_s=2$, using a Gauss-Legendre integrator.
The problem is generated using the Step reformulation and cross-complementarity mode 7 (cf. \cite{Nurkanovic2022c}), and is lifted into the vertical form.

\begin{landscape}
	\begin{longtable}{l l p{6cm} l l p{1.8cm} l}
		\caption{\label{tab:NOSBENCH-ODEs}Continuous time problems used to generate the NOSBENCH test set}\\
		\toprule
		Nr.&Problem Name&Short Description&Type&Class& {Number~of} MPCCs&Citation\\\midrule\endhead
		1&3CPCLS&Two dimensional representation of three carts with only one actuated. &OCP&CLS direct&  18&\cite{Lin2022}\\
		2&3CPTF&One dimensional representation of three carts with only one actuated. &OCP&CLS time-freezing &  18&\cite{Lin2022}\\
		3&CARHYS&Turbo car example with hysteresis. &OCP& HYS time-freezing &  12&\cite{Nurkanovic2022a}\\
		4&CARTIM&Turbo car example with velocity-dependent switch. &OCP&Fillipov System&  24&\cite{Nurkanovic2022}\\
		5&CPWF&Inverted pendulum on a cart with Couloumb friction. &OCP&Fillipov System& 48&\cite{Howell2022a}\\
		6&DAOBCLS&Disc control example.&OCP&CLS direct&  9&\cite{Nurkanovic2023b}\\
		7&DISCM&MDisc control example. &OCP&CLS time-freezing&  6 &\cite{Nurkanovic2023b}\\
		8&DRNLND&Two dimensional control of a landing drone. &OCP&CLS time-freezing&  18 &\cite{Kong2021}\\
		9&DSCOB&Disc control example with a central obstacle, modeled as a nonlinear constraint.&OCP&CLS time-freezing& 27 &\cite{Nurkanovic2023b}\\
		10&DSCSP&Disc control example with discs switching positions. &OCP&CLS time-freezing&  9 &\cite{Nurkanovic2023b}\\
		11&HOPOCP&Hopping robot actuated with a linear leg and reaction wheel in two dimensions. &OCP&CLS time-freezing&  12 &\cite{Howell2022a}\\
		12&MFTOPT&Time optimal control of a linear voice-coil motor. &OCP&Fillipov System&  24&\cite{Christiansen2008}\\
		13&MNPED&Control of a planar monoped robot. &OCP&CLS time-freezing&  36& \cite{Carius2020}\\
		14&MWFOCP&Optimal control of a linear voice-coil motor.&OCP&Fillipov System& 24&\cite{Christiansen2008}\\
		15&SCHUMI&Simple two-dimensional car model with track constraints. Both time optimal and non-time optimal forms.&OCP&Fillipov System&  72&\cite{Stewart2010}\\
		16&SMOCP&Optimal control problem with multiple sliding modes.&OCP&Fillipov System&  48&\cite{Nurkanovic2022}\\
		17&TFBIB&Control of a planar particle in a box with a rotating reference&OCP&CLS time-freezing &  12& \cite{Nurkanovic2023a}\\
		18&TNKCSC&Cascade of tanks with state-dependent switches in flow rate.&OCP&Fillipov System &  6& \cite{Baumrucker2009}\\
		19&TIMF1D&One dimensional particle with contacts.&Simulation&CLS time-freezing&  12&\cite{Nurkanovic2023b}\\
		20&2BCLS&Two balls connected with a stiff spring, that undergo contact dynamics. &Simulation&CLS direct&  9&\cite{Bruls2014}\\
		21&986EQ&Simplified model of structural pounding used in the study of the effects of earthquakes. &Simulation&Fillipov System&  12&\cite{Calvo2016}\\
		22&986FO&Two masses linked by a spring and moving on a surface with Coulomb friction.&Simulation&Fillipov System&  18&\cite{Calvo2016}\\
		23&986FV&Oscillator with varying friction force.&Simulation&Fillipov System&  18&\cite{Calvo2016}\\
		24&986OM&Oscillating mass with contact surfaces.&Simulation&Fillipov System&  12&\cite{Calvo2016}\\
		25&CLS1D&One dimensional particle with contacts.&Simulation&CLS direct&  12&\cite{Nurkanovic2023b}\\
		26&FBS1S& Spring-connected blocks with Coulomb friction. &Simulation&Fillipov System&  18&\cite{Stewart1996}\\
		27&OSCIL&An unstable oscillator with a state dependent switch.&Simulation&Fillipov System&  12&\cite{Nurkanovic2022}\\
		28&RFB1S&Relay feedback system. &Simulation&Fillipov System& 18&\cite{Piiroinen2008}\\
		29&SMCRS&Simple Filippov system with a boundary crossing. &Simulation&Fillipov System&  6&\cite{Nurkanovic2022}\\
		30&SMLSM&Simple Filippov system with state entering a sliding mode.&Simulation&Fillipov System&  6&\cite{Nurkanovic2022}\\
		31&SMSLM&Simple Filippov system with state leaving a sliding mode.&Simulation&Fillipov System&  6&\cite{Nurkanovic2022}\\
		32&SMSPS&Simple Filippov system with state spontaneously leaving an unstable sliding mode.&Simulation&Fillipov System& 6&\cite{Nurkanovic2022}\\
		33&TFPOB&A collapsing pile of balls with contacts&Simulation &CLS time-freezing &  12& \cite{Nurkanovic2022b,Nurkanovic2023a}
	\end{longtable}
\end{landscape}

\subsubsection{Problem subsets}
It takes a lot of CPU time to run the full suite and doing so is not necessary to gain some first insight into the comparative performance of different solution methods.
Therefore, we provide several smaller subsets of problems that can be used to benchmark any future solvers and are used in Section \ref{sec:results} to evaluate existing methods.
\paragraph{Simple problem benchmark - \nosbenchs}
The first subset of \nosnoc\ is a benchmark that only uses {\NS} MPCCs which come exclusively from Fillipov systems and only contain the simplest time-freezing and CLS examples.
These tend to produce relatively easier-to-solve MPCCs and as such this set is an effective way to identify particularly poor-performing solvers.
It contains approximately equivalent numbers of simulation and optimal control problems and typically all of the problems can be solved with the existing state-of-the-art in less than an hour per problem.

\paragraph{Representative small benchmark - \nosbenchrs}
This subset of \nosbench\ is an even smaller but more representative benchmark.
It contains {\NRS} MPCCs that include ones from FESD-J and time-freezing reformulations. 
This subset is primarily used for a second preliminary screening of solvers as it provides more insight into the performance of solvers on problems ranging from the easiest to the most difficult within \nosbench.
\paragraph{Representative large benchmark - \nosbenchrl}
This subset is a set of {\NRL} problems that is made up of a representative sample of all problem difficulties.
This benchmark is meant to be the main problem set that is used to benchmark solvers, and will continue to be expanded as new and interesting problems are added to~\nosbench.
\paragraph{The full benchmark collection - \nosbenchf}
This set contains all 603 MPCCs generated in this version of {\nosbench}.
We aim to test the most competitive solver-method combinations on this set to get a clear performance picture.

\subsection{Stopping criteria and solution quality}
Beyond the NLP stopping criteria used in the underlying NLP solver, we use a further stopping criterion on the homotopy iteration that is based on the ``complementarity residual'' which is defined as:
\begin{align*}
	r_{\perp}(w) = \max( \{ G_i(w)H_i(w) \mid \ \text{for all}\ i \in \{1,\ldots,m\}\}).	
\end{align*}
The stopping condition used is the achievement of a successful NLP solve in the homotopy loop, whose result has a complementarity residual smaller than a given tolerance.
For all the following experiments (except if otherwise noted) we use a complementarity residual tolerance of $\texttt{comp\_tol}= 10^{-7}$.
In the case of \ipopt\ we treat any solution that is reported as optimal or ``solved to an acceptable level'' as a successful solve.
We further accept solutions where the search direction becomes too small if they meet the complementarity tolerance as well as are primal feasible.
We set the runtime limit for a single MPCC solve for all solvers to one hour (3600 seconds) cumulative wall time.

When analyzing the results in the next section, we include in the analysis the quality of a given solution if the problem comes from a discretized OCP.
We do not apply this check for simulation problems as the solutions of such problems are usually unique or at least locally isolated \cite{Nurkanovic2022}.
This verification is done by checking the relative objective value of a given problem-solver pair against the best-known found solution.
We then treat solutions that exceed the best-known objective by at least a factor of two as failures.
This approach is used in order to better evaluate the solution methods specifically in an optimal control context as in this context a significantly worse solution is often a sign of a failure of the solver to achieve the goals of the controller.

\section{Computational results}\label{sec:results}
In this section, several experiments are carried out to evaluate different solution methods for MPCCs.
We first compare the different regularization methods discussed in Section \ref{sec:relaxation_methods} paired with the different homotopy parameter steering strategies from Section \ref{sec:steering} (and Table~\ref{tab:parameter_steering}).
We then explore the space of the homotopy parameters which are used to drive the complementarity regularizations toward the exact complementarity set.
This is followed by a comparison of three NLP solvers (\ipopt~\cite{Waechter2006}, \snopt~\cite{Gill2005}, and \worhp~\cite{Bueskens2013}) used to solve the regularized NLPs.
Moreover, we identify the type of stationary points to which the solver converges. 
In case they are not S-stationary, we solve an LPCC to check if we have a B-stationary point.
Further numerical results on the \nosbench\ test set can be found in the master thesis~\cite{Pozharskiy2023}.
\subsection{General experiment setups}\label{sec:experimental_setup}
All benchmarks are run using an Intel Xeon W-2225 4-core processor with a base clock of 4.1 GHz and a boost clock of 4.6 GHz.
In all cases where we are not explicitly varying the NLP solver, we used \ipopt~\cite{Waechter2006} as the default option.
It is used with its default options except for the following changes: 
\texttt{bound\_relax\_factor = 0}, 
\texttt{mu\_strategy = adaptive}, 
\texttt{mu\_oracle = quality-function},
\texttt{acceptable\_tol =} $10^{-6}$,
\texttt{tol =}$10^{-12}$,
\texttt{dual\_inf\_tol = }$10^{-12}$,
\texttt{comp\_inf\_tol = }$10^{-12}$.
We also default to using \texttt{MA27}~\cite{Duff1982} as the linear solver in \ipopt\ as it has, in our experience, been the most stable of the \texttt{HSL} solvers~\cite{HSL}.
This is due to our observation that both \texttt{MA57}, and \texttt{MA97}, the solvers recommended as state of the art, occasionally cause \ipopt\ crashes due to segmentation faults.
The single-threaded nature of \texttt{MA27} also allows us to run multiple \texttt{IPOPT} instances in parallel in order to improve the throughput of the benchmark.

We measure the performance of each solver primarily using a wall time timer which sums the real time taken to solve all NLPs, and ignore any processing time in between, both as it is not relevant to solver performance, and as it is generally equivalent between different solution methods as it mostly depends on the size of the problem.
The benchmark results are given in the Dolan-Mor\'e performance profiles~\cite{Dolan2002}.
\begin{figure}[t]
	\centering
	\subfloat[Direct NLP solve.]{\includegraphics[width=0.5\textwidth]{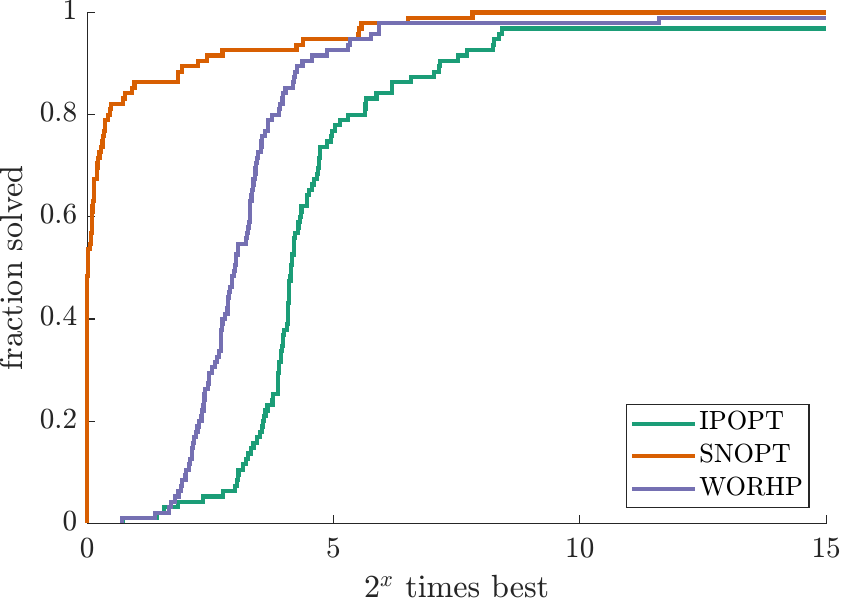}}
	\subfloat[Standard homotopy.]{\includegraphics[width=0.5\textwidth]{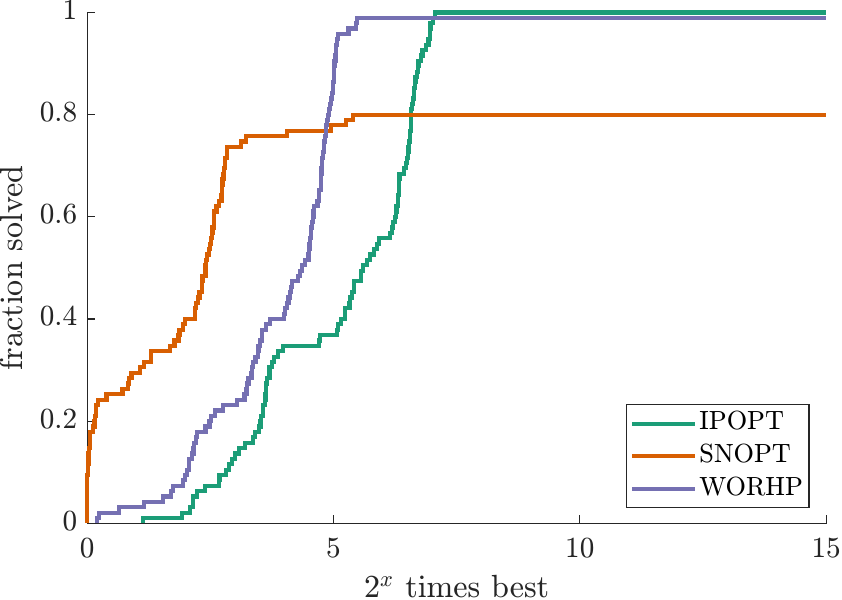}}\\
	\subfloat[$\ell_1$-mode relaxation.]{\includegraphics[width=0.5\textwidth]{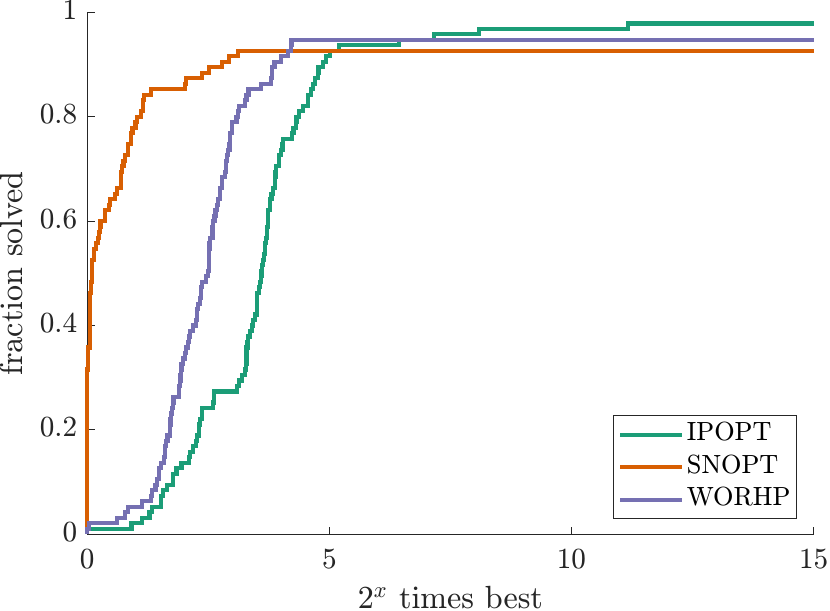}}
	\subfloat[$\ell_\infty$-mode relaxation.]{\includegraphics[width=0.5\textwidth]{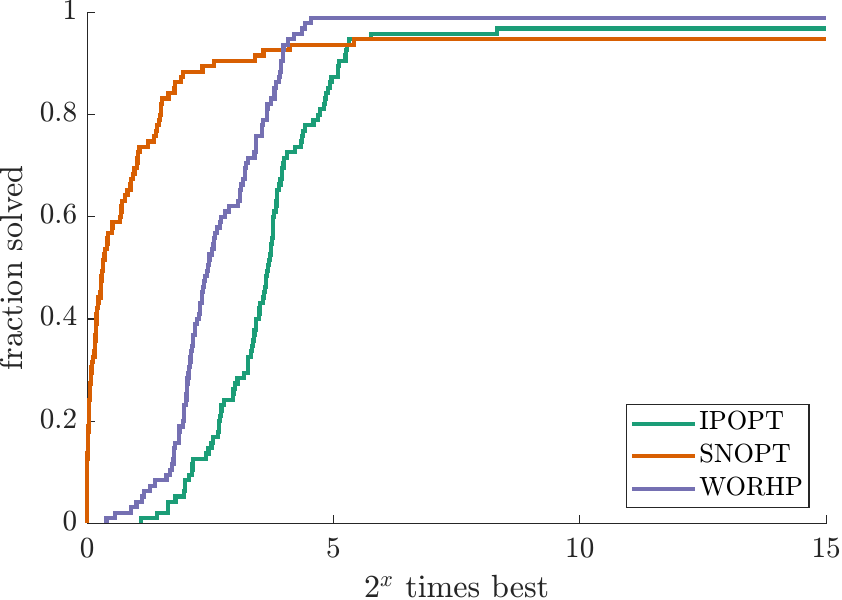}}
	\caption{Evaluating MPCC methods with different NLP solvers on the \MacMPEC\ test set. Each solver-method variant is compared to the other, but the results are split into four plots for better readability. }
	\label{fig:macmpec}
\end{figure}
\subsection{Validating our implementations on MacMPEC}\label{sec:validation}
Before we test the MPCC methods on the \nosbench~collection, we verify the correctness of our the implementation of the regarded MPCC method-NLP solver combinations on the \MacMPEC\ problem set. 
The {\MacMPEC} problem set is available in the form of a tarball of AMPL \cite{Fourer2003} format \texttt{.mod} and \texttt{.dat} files from \url{https://wiki.mcs.anl.gov/leyffer/index.php/MacMPEC}.
We use a modified version of \casadi~\cite{Andersson2019} and {manage} to successfully extract 95 out 106 MPCC problems from this benchmark set.
This is done by first generating \texttt{.nl} files for each problem, then reading these in and generating \casadi~MPCCs of the form in \eqref{eq:mpcc_nosbench}.
We drop any problems that contain complementarities with a ``body'' parameter of 3, as described in \cite{Gay2005}, which are slightly more generic than our implementation permits.
This test suite is then run on four different approaches: 
the direct method (Sec. \ref{sec:mpcc_direct}),
standard Scholtes relaxation  (Sec.\ref{sec:scholtes}), 
$\ell_1$-penalty (Sec. \ref{sec:mpcc_penalty}), 
and the $\ell_\infty$-mode Scholtes relaxation  (Sec.\ref{sec:steering}), 
\color{black}
using the three NLP solvers: \ipopt~\cite{Waechter2006}, \snopt~\cite{Gill2005}, and \worhp~\cite{Bueskens2013}.
These were the overall most successful variants also in the later experiments, and for {the} sake of brevity we do not run the benchmark on all possible solver-method combinations.

We note that several of these approaches solve each of the 95 problems in under ten minutes, and most solve more than 90\% of the problem set in the same amount of time.
These results are summarized in Figure~\ref{fig:macmpec}. 
They align with the results reported in other {papers}~\cite{Thierry2020,Izmailov2012a,Fletcher2002a}, which validates the correctness of our implementation.

In contrast to what we will see in the subsequent sections, we observe better performance of the direct method on the smaller problems in this test set, along with much better performance from \snopt, again tied to the smaller size of the problems.
In general this supports our assertions on the relative difficulty of the \nosbench\ test suite when compared to existing state-of-the-art benchmarks.

\subsection{Evaluating different MPCC methods}\label{sec:mpcc_methods_benchmark}
In this section, we compare nine variants of the relaxation-based methods discussed in Section \ref{sec:relaxation_methods}: Scholtes' relaxation, three smoothed NCP functions (Fischer-Burmeister (FB), Natural-Residual(NR) and Chen-Chen-Kanzow (CCK), cf. Section \ref{sec:scholtes}), the Lin-Fukushima (LF) method.
We also test the Steffensen-Ulbrich method with the two test functions mentioned in Section \ref{sec:steffensen_ulbrich} (denoted by SU1 and SU2), and the kinked relaxations by Kadrani et. al and Kanzow-Schwartz (KS).
Each of the nine methods is tested along with one of the three methods for steering the relaxation parameter summarized in Table~\ref{tab:parameter_steering}, which results in a total of 27  different methods.
Furthermore, we solve the MPCCs directly as NLPs (cf. Section \ref{sec:mpcc_direct}) and with an exact-$\ell_1$-penalty method without slacks (cf. Eq. \eqref{eq:mpcc_ell_1_penalty}).
An implementation of all aforementioned methods is available in \nosnoc.

As we discussed in Section \ref{sec:steering}, the algorithms implemented in \nosnoc\ for solving MPCCs have several free parameters, namely $\sigma_0$, $\kappa$, and $\eta$.
In this experiments we use \eqref{eq:parameter_homotopy} for the $\sigma_k$ updated with $\sigma_0 = 1, \kappa = 0.1$. 
In the subsequent sections, we vary these parameters in order to evaluate a good set of default parameters.
The comparison is run on the \nosnoc-\texttt{S} subset.

\begin{figure}[t]
	\centering
	\subfloat[Standard relaxation.]{\includegraphics[width=0.5\textwidth]{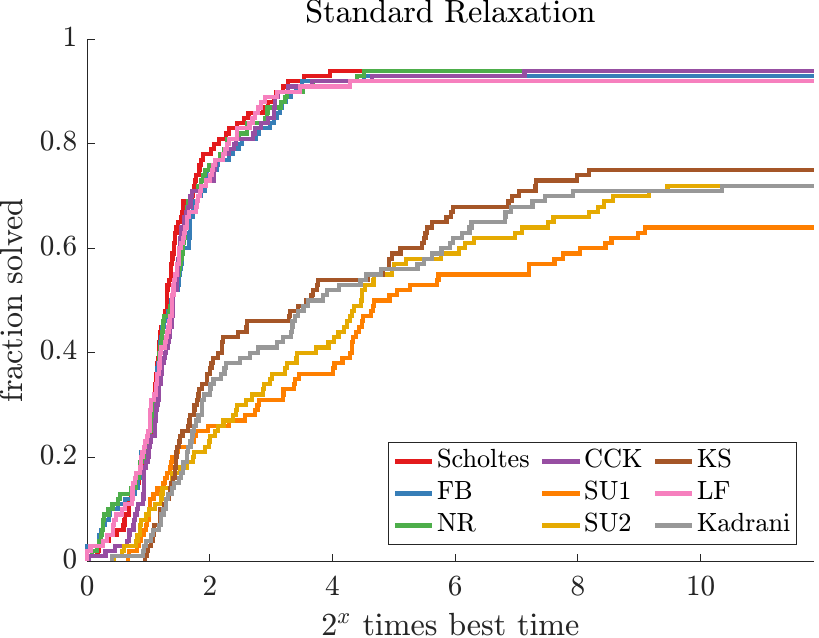}}
	\subfloat[$\ell_\infty$-mode relaxation.]{\includegraphics[width=0.5\textwidth]{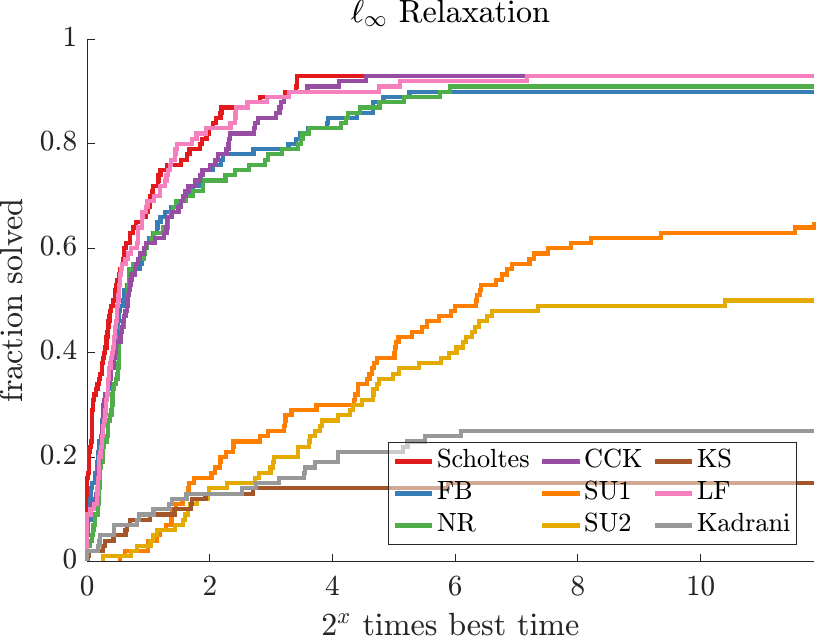}}\\
	\subfloat[$\ell_1$-mode relaxation.]{\includegraphics[width=0.5\textwidth]{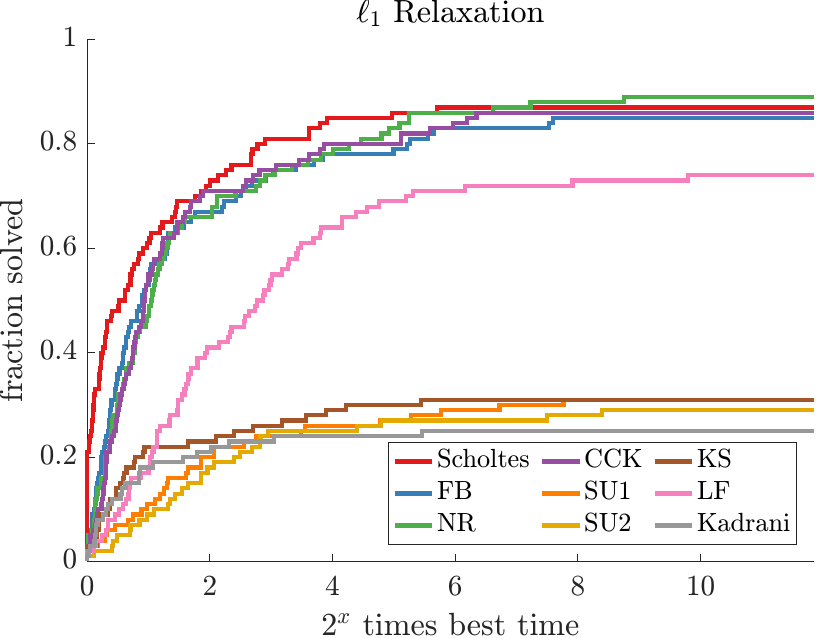}}
	\subfloat[Direct method and $\ell_1$-penalty.]{\includegraphics[width=0.5\textwidth]{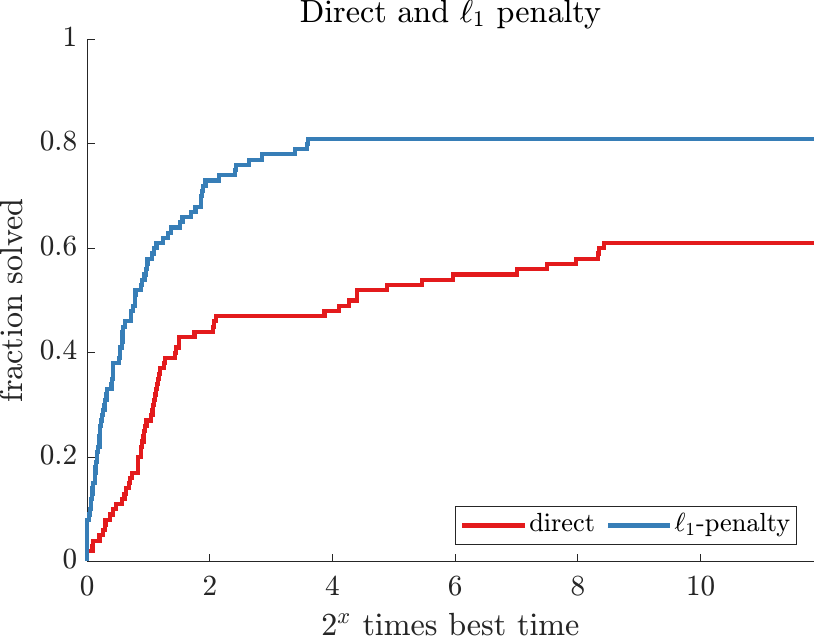}}
	\caption{Evaluating regularization methods for MPCCs on the \nosbench-\texttt{S} subset.}
	\label{fig:regularization_profile}
\end{figure}
\begin{figure}[t]
	\centering
	\includegraphics[width=\textwidth]{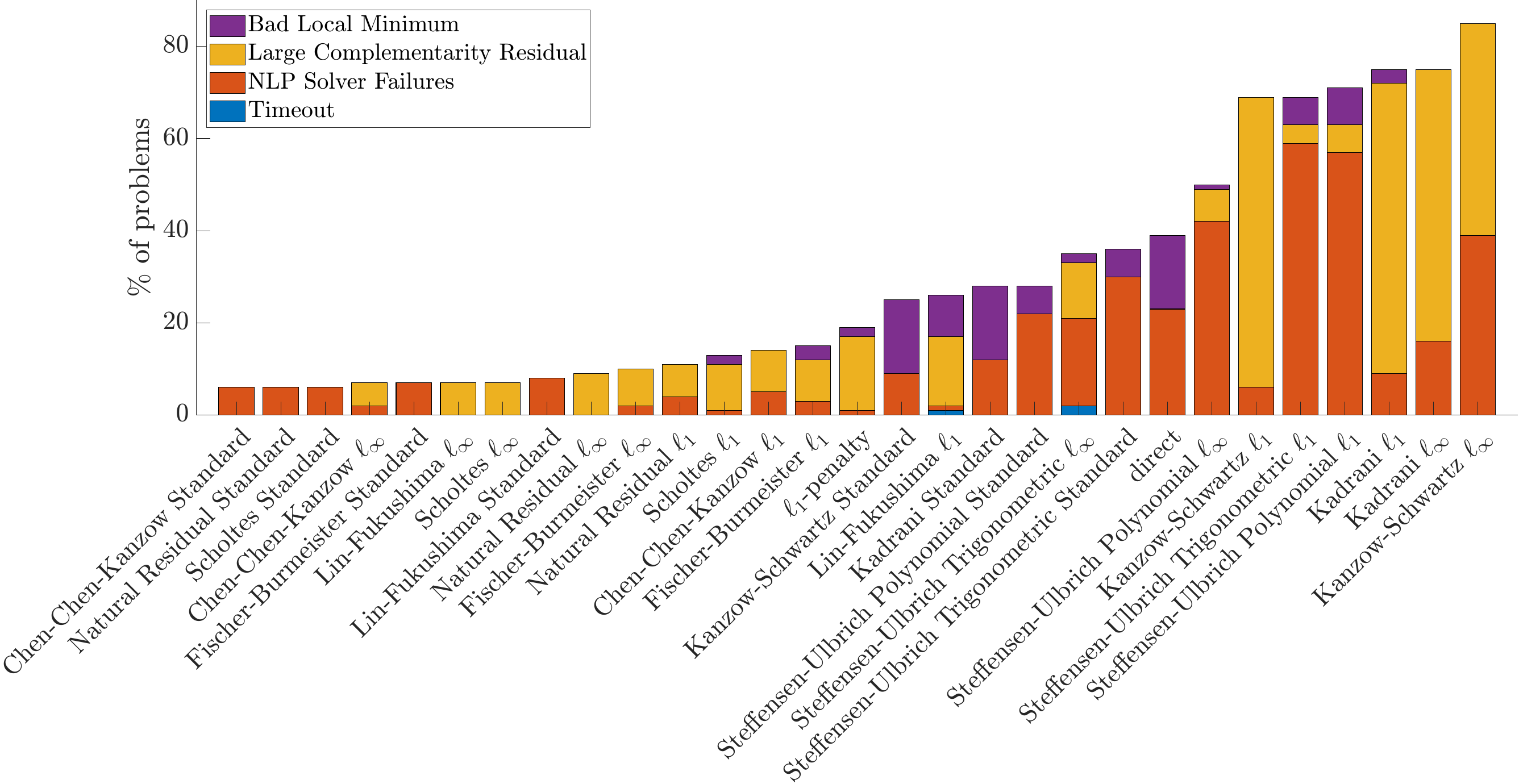}
	\caption{Failure reasons for different regularization methods on the \nosbench-\texttt{S} subset.}
	\label{fig:regularization_failure}
\end{figure}

An overview of the relative performance of each MPCC solution method with respect to all others can be seen in Figure \ref{fig:regularization_profile}.
The comparison is split into four subplots for better readability.
Figure \ref{fig:regularization_failure} summarizes the reasons for the failure of each of the 29 methods we test.

The general outcome of this benchmark is a victory for the Scholtes relaxation and the smoothed NCP functions (which lead to the same feasible set as Scholtes' method) in all of its three steering strategies.
From the performance plots, we can clearly see that for all three methods of steering the relaxation parameter, the Scholtes relaxation successfully solves almost as many or more problems than the other relaxation methods.
On the other hand, the more sophisticated local and nonsmooth relaxation methods perform worse in all cases. 
This complies with the fact that these methods have weaker theoretical properties \cite{Kanzow2015} if the subproblems are solved inexactly (which is inevitable in finite precision arithmetic).

We can in detail compare the relaxations using the standard approach to drive the regularization parameter to zero.
The results are depicted in Figure \ref{fig:regularization_profile} (a).
Here we see a performance lead for the Scholtes relaxation, albeit a very slim one.
Methods arising from the standard parameter steering can be split into three different groups based on performance (in order): the group containing the Scholtes relaxation as well as those that use NCP functions, Lin-Fukushima and Kanzow-Schwartz which are nearly as fast but plateau earlier (i.e., solve fewer problems),
and the Steffensen-Ulbrich and Kadrani relaxations which gain only about 10\% robustness over the direct method but are somewhat slower.
We note that Figure \ref{fig:regularization_profile} suggests that we see very few outright victories for the standard Scholtes method, however, it is the most robust, achieving the highest fraction of problems solved, 95\%.

We then move on to our analysis of the $\ell_\infty$-mode relaxation with the results given in Figure \ref{fig:regularization_profile} (b).
One can see that the Scholtes method paired with this parameter steering strategy wins, with the largest fraction of successful solves.
It maintains its lead but only reaches a solution on about 92\% of the problems.
Once again we see that the NCP function relaxations approach the performance of the Scholtes relaxation and match it in robustness.
We also see the Lin-Fukushima relaxation perform well again but not quite at the level of the Scholtes-type group.
On the other hand, we see extremely poor performance from the Steffensen-Ulbrich relaxations and the kinky relaxations of Kadrani et al. and Kanzow-Schwartz.
From Figure \ref{fig:regularization_failure} we see that both the Kadrani and Kanzow-Schwartz relaxations fail frequently with a point at which the NLP solver claims it is optimal, but still has a large complementarity residual.

Finally, we analyze the performance of the $\ell_1$-mode relaxations and the $\ell_1$-penalty formulation, given in Figures \ref{fig:regularization_profile} (c) and (d), respectively.
The $\ell_1$-penalty method nearly ties the $\ell_\infty$ Scholtes relaxation in absolute wins, however, it fails to break the 90\% of problems solved.
This performance is nearly mirrored by the $\ell_1$-mode Scholtes relaxation albeit without the outright victories of the penalty method.
Here we also see the second major gap between the Scholtes relaxation and the NCP function-based relaxations, though they remain the best-performing group.
Surprisingly, the remaining relaxations perform worse than even the direct NLP solve approach.

Another notable result is that while for about half the problems in this set the direct NLP approach solves the problem in an acceptable time frame, it quickly stalls at that point, cf. Figure \ref{fig:regularization_profile} (d).
Its failures, as seen in Figure~\ref{fig:regularization_failure}, are fairly evenly split between converging to unacceptable local minima, and failing to converge at all, primarily due to either step calculation failures at unacceptable points, or due to claiming infeasibility.
This already points to the usefulness of the homotopy methods as it is clear to see that those methods perform better on a large proportion of problems with little to no impact on absolute performance.

\begin{figure}[b]
	\centering
	\subfloat[Standard relaxation.]{\includegraphics[width=0.5\textwidth]{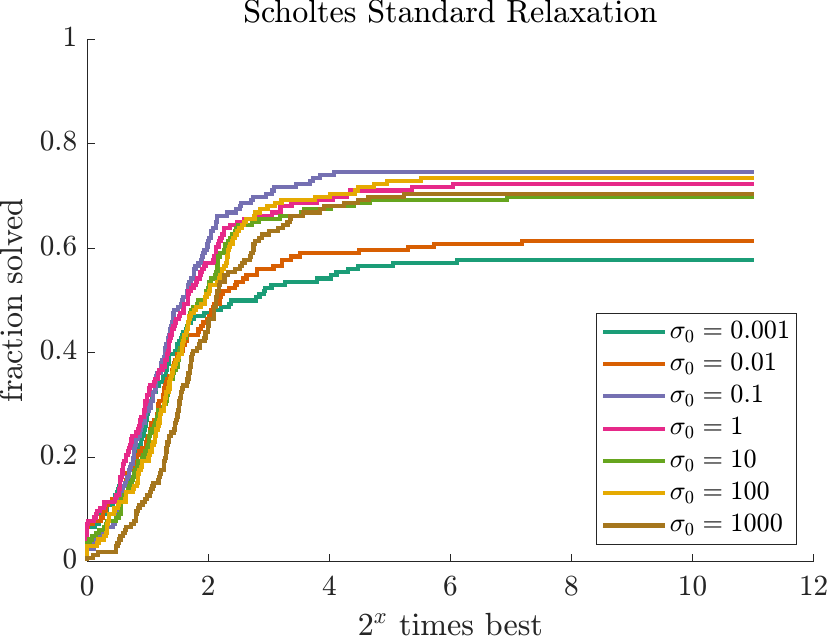}}
	\subfloat[$\ell_\infty$-mode relaxation.]{\includegraphics[width=0.5\textwidth]{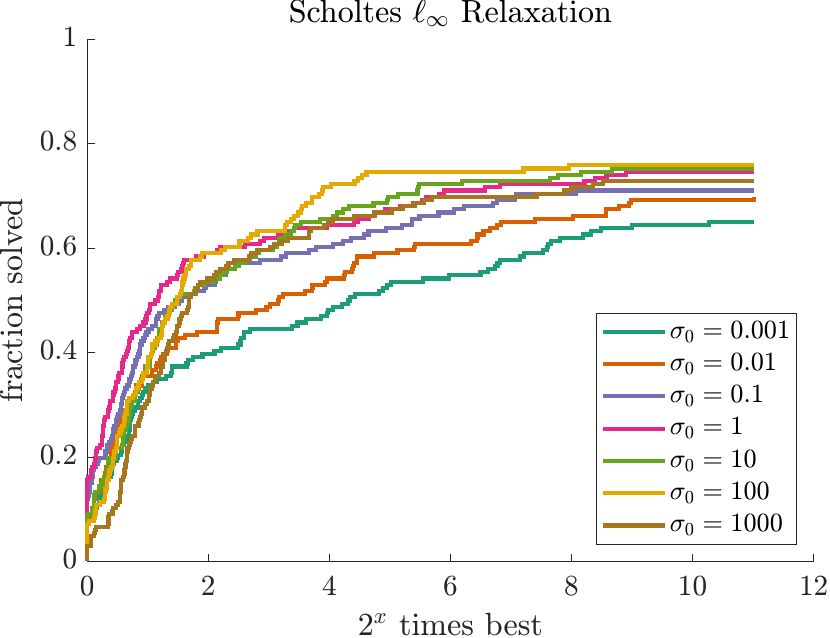}}
	\caption{Evaluating the influence of $\sigma_0$, the initial homotopy parameter in the \nosbenchrl\ test set.}
	\label{fig:sigma_0_experiment}
	\centering
	\includegraphics[width=0.95\textwidth]{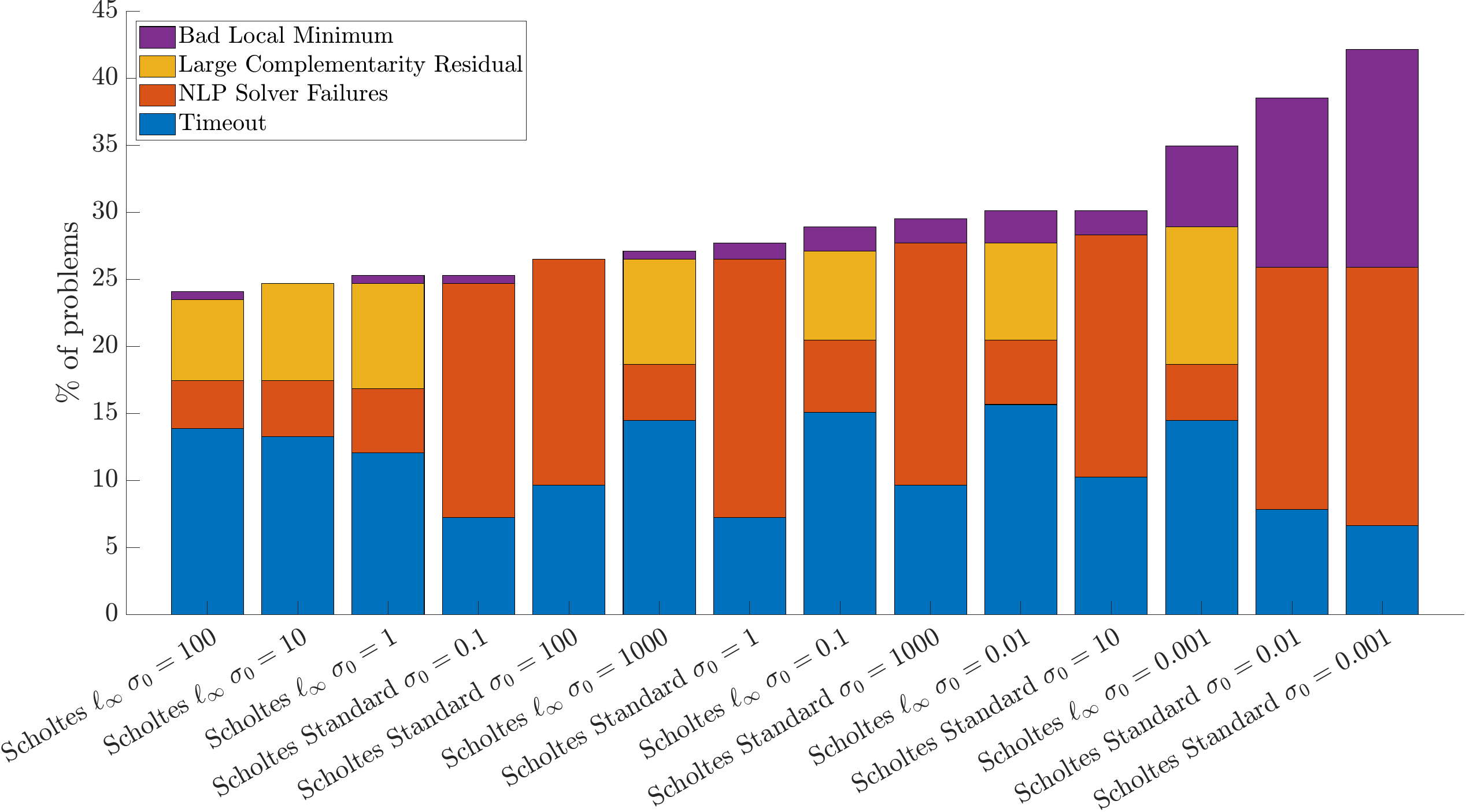}
	\caption{Failure reasons for different values of $\sigma_0$.}
	\label{fig:sigma_0_failure_reasons}
\end{figure}

We conclude that the best methods in terms of speed and robustness are the $\ell_\infty$-mode Scholtes relaxation and the Scholtes relaxation with standard homotopy parameter steering, respectively.
The latter is certainly the choice for robustness as we will show in future experiments. 
The robustness can even be improved by adjusting the $\sigma_0$ and $\kappa$ parameters which govern the trajectory of the regularization parameter.

\subsection{The role of the homotopy parameters}
As described in the previous section, we note that the Scholtes relaxation is the optimal choice in its standard and $\ell_\infty$ forms.
We do not examine further the $\ell_1$ strategy, as it is not competitive with the other two.
In order to elucidate the role of homotopy parameter update strategies, we run an experiment on \nosbenchrl\ (representative large subset), first varying the initial regularization parameter $\sigma_0$ and then varying $\kappa$, the rate at which we drive the regularization parameter to zero.
It turns out that these two parameters have a moderate influence on the stability and speed of the homotopy solver converging to an acceptable solution.

We first discuss the effect of the initial value of the homotopy parameter $\sigma_0$ for which the performance plots can be seen in Figure~\ref{fig:sigma_0_experiment}.
Once again, we compare each method with the other, but split it into two subplots for readability.
The first major takeaway from the analysis of the performance plots is that this has a much less significant effect on the $\ell_\infty$-mode relaxation.
Intuitively this does make sense as in this mode we simply use a penalty factor $\frac{1}{\sigma_k}$ in order to drive the complementarity residual to zero.
Thus, $\sigma_k$ is not a limiting factor on the complementarity residual in each step, depending of course on the scaling of the problem at hand.
However, minimal effect is not no effect and we still see a worsened convergence if we choose a $\sigma_0$ that is too small.
In contrast, it is observed that for the standard relaxation, for very small $\sigma_0$ the solvers converge to points with a complementarity residual exactly equal to $\sigma$.
It is notable that this reduction in performance primarily comes from convergence to significantly worse local minimizers as seen in Figure~\ref{fig:sigma_0_failure_reasons}.

On the other hand, we see much earlier and much more pronounced decay in performance for the standard Scholtes regularization.
We see almost a 30\% reduction in the number of problems solved if $\sigma_0$ is chosen to small. 
Moreover, from Figure \ref{fig:sigma_0_failure_reasons} we see that the primary reason for failure is the NLP solver converging to bad local minimizers, compared to the more successful cases, e.g. with $\sigma_0$ equal to 10 or 0.1.
One possible reason is that, for larger values of $\sigma_0$ the problems are more relaxed and the solver may get easier attracted by better minimizers.
In conclusion, for $\sigma_0$ in the rage of $0.1$ to $10$ the best performance is achieved.

\begin{figure}[t]
	\centering
	\subfloat[Standard relaxation.]{\includegraphics[width=0.5\textwidth]{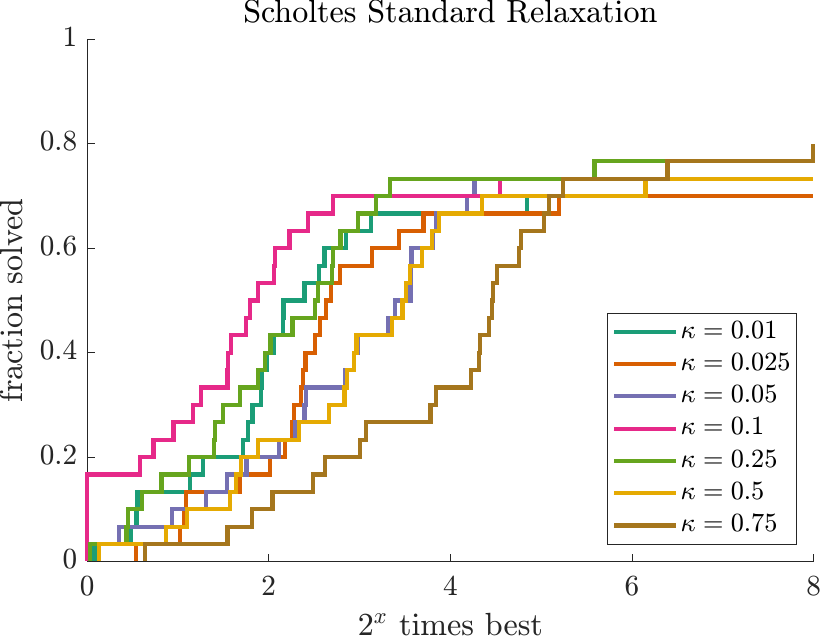}}
	\subfloat[$\ell_\infty$-mode relaxation.]{\includegraphics[width=0.5\textwidth]{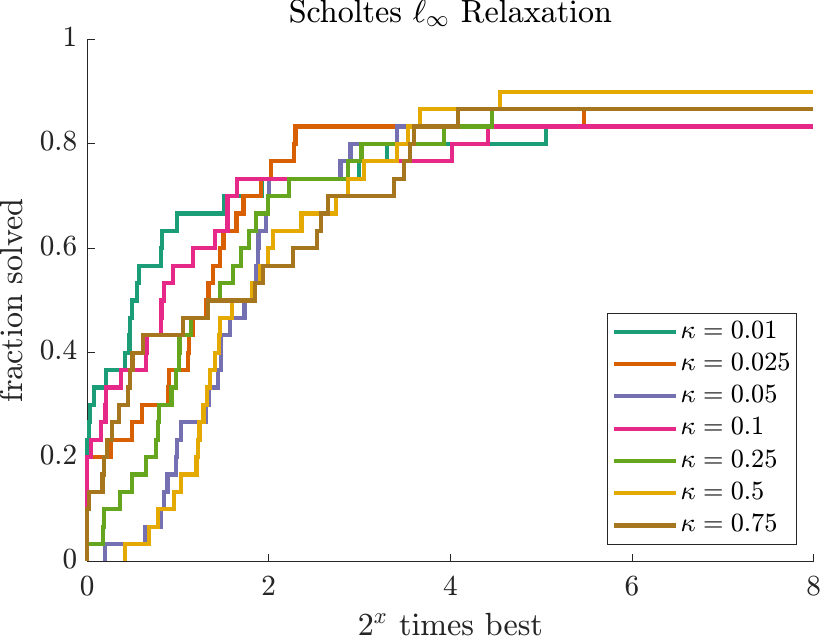}}
	\caption{Evaluating the influence of the homotopy slope parameter $\kappa$ in the \nosbenchrl\ test set.
	}
	\label{fig:sigma_slope_experiment}
\end{figure}

Next, we examine the effect of the homotopy parameter update factor $\kappa$. 
In the literature, several different choices are used. 
Examples are $\kappa = 0.2$ in \cite{Kadrani2009}, $\kappa = 0.01$ in \cite{Scholtes2001} and \cite{Hoheisel2013}, and $\kappa = 0.1$ in \cite{Steffensen2010}.
In this experiment, we fix $\sigma_0 = 1$.
In Figure \ref{fig:sigma_slope_experiment} one can see that the effect of $\kappa$ on the overall success is surprisingly small compared to the initial regularization parameter $\sigma_0$.
In particular, we see essentially no difference in performance for the $\ell_\infty$-mode regularization.
This is very likely due to the fact that for a majority of problems, we see only several (and occasionally only one) homotopy iterations before the solver converges to an acceptable complementarity residual.
The standard relaxation clearly shows both the weakness and the strength of a relatively slow homotopy.
We note that for smaller $\kappa$ the problem converges quicker due to having to take less homotopy iterations to converge to a sufficiently accurate solution.
We also see that this benefit disappears as we get to more difficult problems.
On the other hand, with larger $\kappa$, we have to solve more problems, but they can be solved quicker since the initial guess of the previous solution is much better.
We also see a mild improvement (around 10\%) in the number of problems that are solved with the larger $\kappa$.
This makes it likely that a larger homotopy update rate is particularly useful for ensuring convergence for more difficult problems while a smaller $\kappa$ (or use of the $\ell_\infty$-mode) is the superior option for simpler problems.
To conclude, a good default choice is $\kappa=0.1$ and $\sigma_0 = 1$ with the $\ell_{\infty}$ mode.

\subsection{Evaluating different NLP solvers}
\begin{figure}[t]
	\centering
	\subfloat[Standard regularization.]{\includegraphics[width=0.5\textwidth]{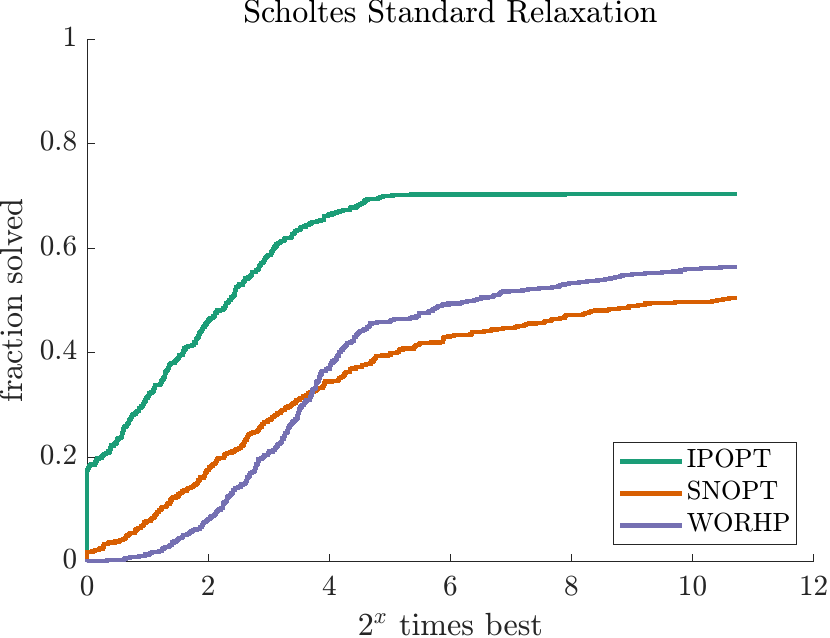}}
	\subfloat[$\ell_\infty$-mode regularization.]{\includegraphics[width=0.5\textwidth]{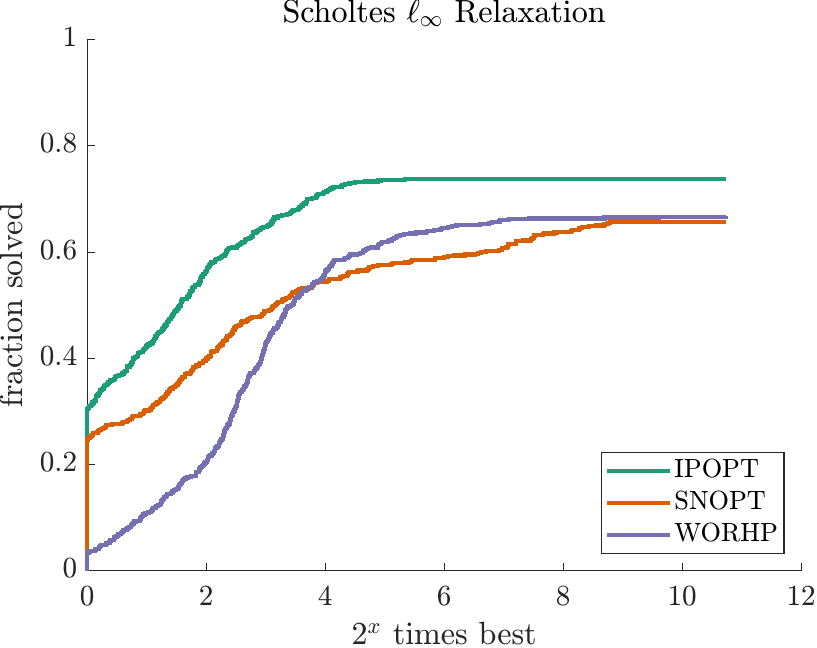}}
	\caption{Evaluating different NLP solvers on the \nosbenchf\ test set.
	}
	\label{fig:nlp_solvers}
	\centering
	\includegraphics[width=0.84\textwidth]{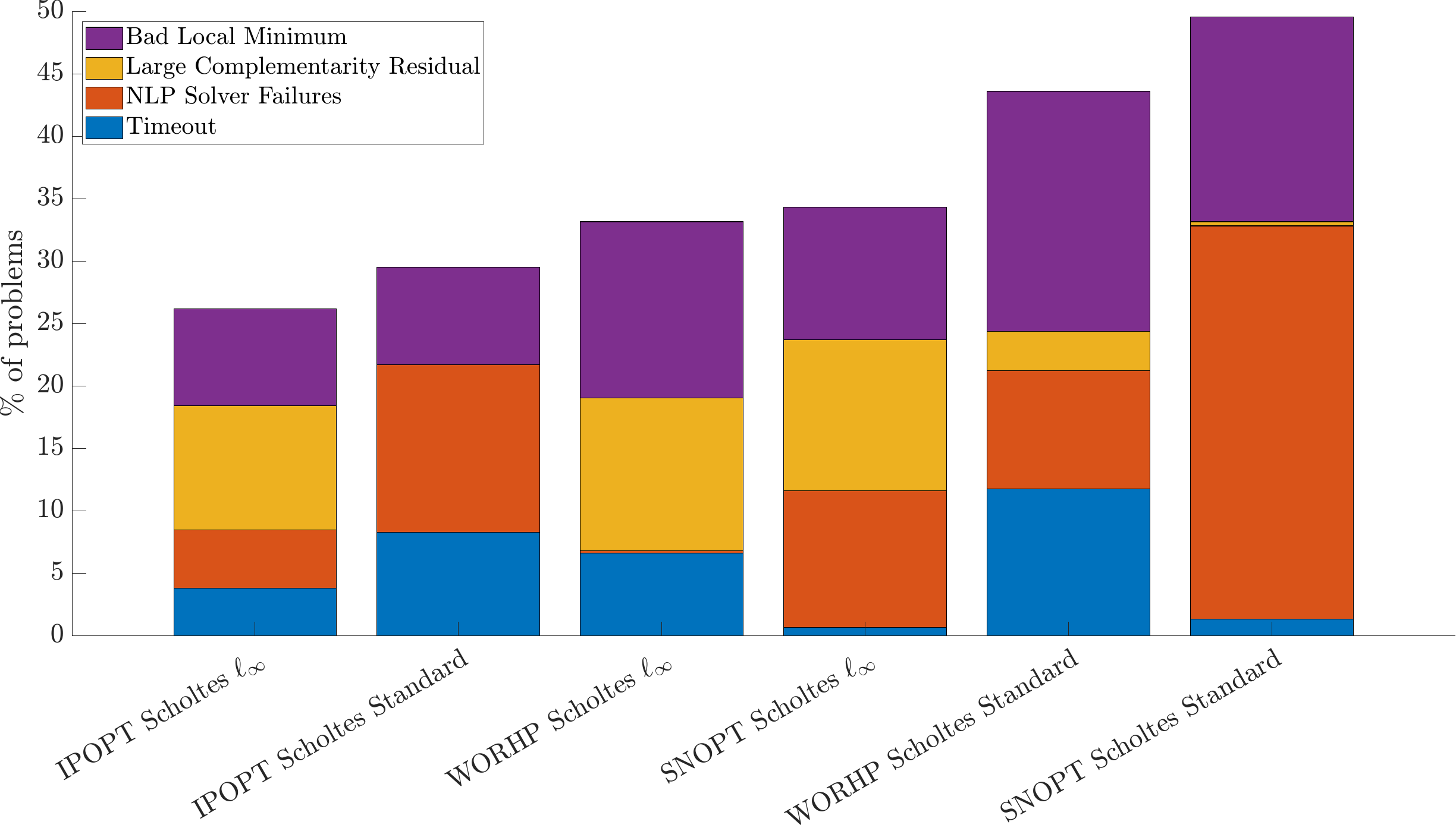}
	\label{fig:nlp_solvers_failures}
	\caption{Failure reasons for different NLP solvers.}
\end{figure}
So far we used in our experiments only \ipopt\ as NLP solver in the homotopy loop. 
Next, we compare this solver to \snopt~\cite{Gill2005}, and \worhp~\cite{Bueskens2013} on the \nosbenchf\ set.

All three solvers are tested on the Scholtes relaxation with the standard and $\ell_{\infty}$ relaxation parameter steering strategies.
In this case, we use the existing default homotopy parameters $\sigma_0 = 1$ and $\kappa=0.1$ and use the standard default settings for both \snopt\ and \worhp, except for those related to maximum iterations and timeouts which are set as for \ipopt.
The solvers are tested on the full \nosbenchf\ test suite, which contains \NF\ MPCCs.
Figure~\ref{fig:nlp_solvers} shows the performance profiles for the three solvers for the two different parameter steering strategies.

On some easier and smaller problems \snopt\ is the fastest solver, which is expected due to the specialization of active set methods on small to medium size problems.
\worhp\ on the other hand is is not particularly fast but is more robust than \snopt, and this gap is more prevalent in the standard relaxation.
It is however clear to see that \ipopt\ is by far the winner here with both the most overall wins in the $\ell_\infty$-mode where it successfully solves 73.8\% of the problem set.
This may be caused by slightly optimized solver settings we have used for \ipopt. 
Moreover, the homotopy parameters $\kappa = 0.1$, $\sigma_0 = 1$ are also somewhat tuned for \ipopt. 
Other values might be beneficial to the performance of the other solvers.

In Figure \ref{fig:nlp_solvers_failures} we also report the reasons for the failure of the NLP solver on the test set. 
Some optimal control problems have very nonlinear dynamics and combined with the relaxed complementarity constraints one obtains difficult NLP subproblems.
Moreover, in the $\ell_{\infty}$ approach, the slack variable and consequently the complementarity residual cannot be brought to a sufficiently small value, despite a very large penalty parameter in the objective.
We note also that some of the problems in \nosbenchf\ are very large, and they might get solved if the solvers were allowed more time.
Interestingly, in more then 5\% of cases \ipopt\ fails because the other approaches have found significantly better local minima. 
This is usually the case for MPCCs coming from OCPs with nonsmooth systems with state jumps.

\paragraph{Type of stationary points}
\begin{figure}[t]
	\centering
	\subfloat[Shares of problems with empty ($|\mathcal{I}_{00}|=0$) and nonempty ($|\mathcal{I}_{00}|>0$) biactive set.]{\includegraphics[width=0.48\textwidth]{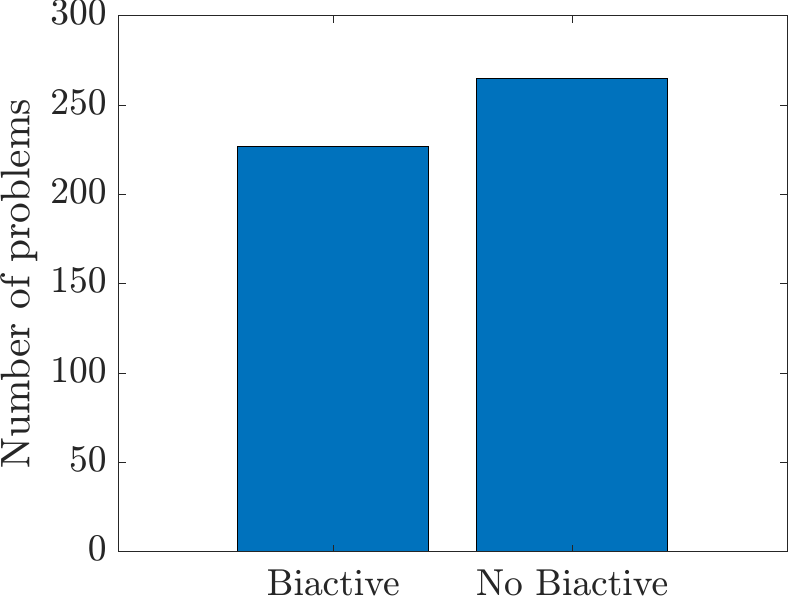}}
	\hfill
	\subfloat[Distribution of stationary points defined in Def.~\ref{def:mpcc_stationarity}, ND stands for not decided.
	]{\includegraphics[width=0.48\textwidth]{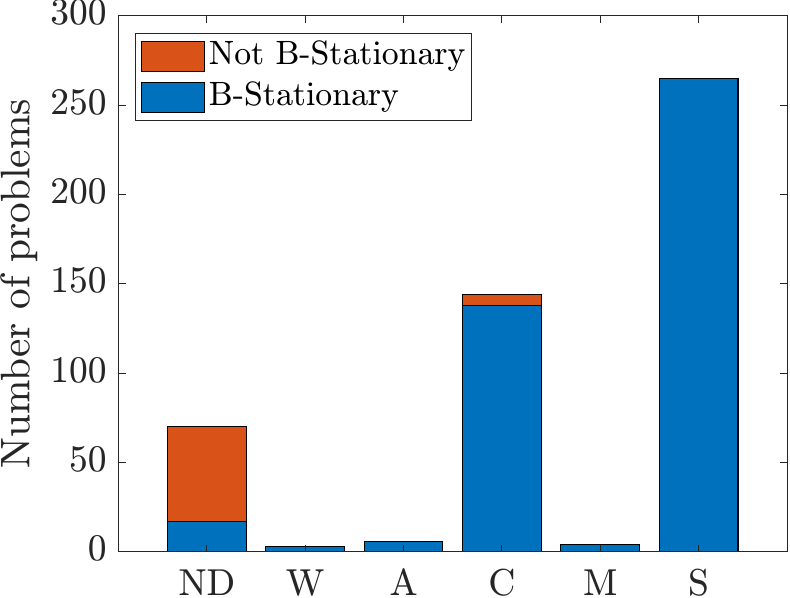}}
	\caption{Evaluating type of stationarity points in the \nosbenchf\ test set.
	}
	\label{fig:stationarity_types}
\end{figure}
We report also the statistics of the type of stationarity points for the most successful solver-method combination, namely \ipopt\ with the Scholtes relaxation and $\ell_{\infty}$-parameter steering.
he set $\mathcal{I}_{00}$ is often called the biactive complementarity constraints
Figure \ref{fig:stationarity_types} 
(a) shows the number of problems with an empty ($|\mathcal{I}_{00}|=0$) and nonempty set ($|\mathcal{I}_{00}|>0$) of biactive complementarity constraints.
Recall that problems with a biactive set are more regular. 
The biactive set is calculated via checking the values of $G(w)$ and $H(w)$. 
We first calculate the biactive set $\mathcal{I}_{00}$ as all $i$ such that $G_i(w)<\sqrt{\texttt{comp\_tol}}$ and $H_i(w)<\sqrt{\texttt{comp\_tol}}$. 
We then calculate $\mathcal{I}_{0+}$ as all $i$ not in $\mathcal{I}_{00}$ that satisfies $G_i(w) < H_i(w)$ and $\mathcal{I}_{+0}$ as all $i$ not in $\mathcal{I}_{00}$ that satisfies $G_i(w) > H_i(w)$.
This however sometimes fails to correctly identify the active set which leads the TNLP to fail to converge. 
In these cases we iteratively remove elements from $\mathcal{I}_{00}$, with the pairs $(G_i(w),H_i(w))$ that are furthest from the origin being removed first and added to the corresponding set based on the magnitude of the components. 
We set the maximum number of iterations in this procedure to $|\mathcal{I}_{00}|$.
In the cases where the TNLP is infeasible or leads to a significantly different objective than obtained by the homotopy approach indicates that we did not identify the active-set correctly. 
We denoted these cases as Not Decided (ND).

Recall that if the biactive set is empty, then the solution is automatically S-stationary.
For those with a nonempty biactive set, we solve the corresponding TNLP (cf. Definition \ref{def:mpcc_nlps}) and asses the type of the stationary point according to Definition \ref{def:mpcc_stationarity}.
The results are summarized in Figure \ref{fig:stationarity_types} (b).
It turns out none of points with a nonempty $\mathcal{I}_{00}$ is S-stationary.
Next, we check additionally for all these points if they are B-stationary, i.e., we check if they permit first-order descent directions.
This can be done by solving the nonconvex LPCC in \eqref{eq:b_stationariry}.
We reformulate the LPCC into an equivalent mixed-integer linear program (MILP)~\cite{Facchinei2003} and solve it with \texttt{Gurobi}~\cite{Gurobi}.
Note that for the verification we need $d = 0$ to solve \eqref{eq:b_stationariry}.
However, this point can be either a local or global solution, but the MILP approach always finds a global minimizer.
To address this, we add a trust region constraint $ \| d\|_{\infty} \leq 10^{-2}$ to make sure that we isolate the (possibly) local optimum $d = 0$. 
If this constraint becomes active, we shrink the trust region radius down to $10^{-4}$, and if $d =0$ is still not optimal, we conclude that it is not a local optimum and that the point $w^*$ is not B-stationary.
Moreover, as sanity check we confirm with the LPCC approach that all S-stationary points are indeed B-stationary.
It turns out that in almost all cases where we could identify the stationary point the points are B-stationary, except in a few cases of C-stationary points.
This provides evidence that MPCCs obtained from nonsmooth optimal control problems, despite sometimes being difficult to solve, are not very degenerate. 
On the other hand, for the ND problems where we could not identify the active set in the TNLP, first-order descent directions often exist. 
One reasons could be that the homotopy loop terminated {too} early. 
In some cases, we noticed that lowering the complementarity tolerance would either make the solver fail, or help to identify the more accurate solution as B-stationary.

\section{Conclusion and outlook}\label{sec:conclusions}
The goal of this paper was to create a benchmark collection of Mathematical Programs with Complementarity Constraints (MPCCs) obtained from the time-discretization of nonsmooth simulation and Optimal Control Problems (OCPs) and to use it for the evaluation of tailored MPCC solution methods.
This provided a large source of practical MPCCs. 
However, MPCCs violate standard constraint qualifications and thus require specialized first-order optimality conditions and solution methods, which were reviewed in detail.
The literature reports very good numerical performance of standard MPCC methods on existing benchmarks. 
Unfortunately, we have not observed such robust performance on MPCCs obtained from nonsmooth OCPs.

To better assess the limitations of the current state-of-the-art and to motivate further development of MPCC methods, we introduce a new benchmark set, which we call \nosbench.
The novel benchmark consists of a total of 603 problems.
Moreover, we derive several subsets from it to facilitate the analysis of a variety of solution methods.
All the methods we test solve a sequence of regularized nonlinear programs (NLPs) in a homotopy approach.
We compare different regularization strategies, different NLP solvers, different approaches to controlling the degree of relaxation in the subproblems, and the influence of the homotopy meta-parameters.
We find that the oldest and simplest methods, namely Scholtes' global relaxation~\cite{Scholtes2001} and smoothed nonlinear complementarity functions~\cite{Facchinei1999} (which are often equivalent to Scholtes' approach), perform best.
This is consistent with previous extensive experiments such as those performed by \cite{Hoheisel2013}.
Surprisingly, our implementations of the more sophisticated regularization strategies show quite disappointing results even on easier subsets of our test set.
Two positive results are that solutions of nonsmooth OCPs are quite often S-stationary and that the weaker stationary points rarely allow first-order descent directions.
However, in the best case, we only manage to solve {\Nsucess} of the problems on the full problem set, which is not yet satisfactory.

It would be interesting to test some of the active-set MPCC methods in the future, should robust open-source implementations of them become available. 
We aim to extend the \nosbench\ test suite by further challenging nonsmooth simulation and optimal control problems.

In conclusion, relaxation-based MPCC methods, coupled with a robust NLP solver, can perform reasonably well even on large and highly nonlinear problems.
The experiments performed in this paper helped us to extract some rules of thumb for the default solver settings in \nosnoc.
However, there is still room for improvement, and the benchmark collection introduced in this paper can help to test novel methods.

\bmhead{Acknowledgments}
This research was supported by DFG via Research Unit FOR 2401, projects 424107692 and 525018088, by BMWK via 03EI4057A and 03EN3054B, and by the EU via ELO-X 953348.


\end{document}